\documentclass[fleqn,reqno,a4paper]{amsart}
\usepackage{style}

\newcommand\place{\mathord{-}}                  
\newcommand\ot{\otimes}                         
\renewcommand\*{\bullet}                        
\newcommand\epi{\twoheadrightarrow}             
\newcommand\mono{\rightarrowtail}               
\newcommand\denotes{\doteq}                     
\newcommand\abs[1]{{\lvert#1\rvert}}            

\newcommand\op{\mathrm{op}}                     
\newcommand\odd{\mathrm{odd}}                   
\newcommand\ev{\mathrm{ev}}                     
\newcommand\id{\mathrm{id}}                     

\newcommand\I{{}'\!}                            
\newcommand\II{{}''\!}                          

\newcommand\pre[3]{\leftidx{_#1^{\vphantom{q}}}{#2}{#3}} 

\newcommand\quot{\mathop{\kern 0pt\mathsf{q}}\nolimits}
\newcommand\rem{\mathop{\kern 0pt\mathsf{r}}\nolimits}

\newcommand\category{\mathsf}                   
\newcommand\Alg{\category{Alg}}                 
\newcommand\Mod{\category{Mod}}                 

\newcommand\lMod[1]{\leftidx{_{#1}}\Mod}                
                          
\newcommand\bMod[2]{\leftidx{_{#1}}\Mod_{#2}}           
\newcommand\Bimod[1]{\leftidx{_{#1}}\Mod_{#1}}
\newcommand\Modf{\category{mod}}
\newcommand\lModf[1]{\leftidx{_{#1}}\Modf}

\newcommand\NN{\mathbb{N}}                      
\newcommand\ZZ{\mathbb{Z}}                      

\let\oldO\O                                     
\renewcommand\O{\mathcal{O}}                    
\newcommand\R{\mathscr{R}}                      
\renewcommand\r{\mathfrak{r}}                   

\newcommand\HH{H\!H}                            

\DeclareMathOperator{\Aut}{Aut}         	
\DeclareMathOperator{\Der}{Der}         	
\DeclareMathOperator{\Ext}{Ext}         	
\DeclareMathOperator{\Tor}{Tor}         	
\DeclareMathOperator{\Op}{Op}           	
\DeclareMathOperator{\sOp}{sOp}         	
\DeclareMathOperator{\im}{im}           	
\DeclareMathOperator{\pdim}{pdim}       	
\DeclareMathOperator{\Z}{Z}             	
\DeclareMathOperator{\rad}{rad}         	
\DeclareMathOperator{\Ch}{Ch}           	

\def\dash---{\thinspace---\hskip.16667em\relax} 
\newcommand\definiendum[1]{\emph{#1}}           
\newcommand\nbd{\nobreakdash}                   
\newcommand\ie{\emph{i.e.}}
\newcommand\cf{\emph{cf.}}
\newcommand\etal{\emph{et al.}}

\begin{document}

\title[The change-of-rings spectral sequence]{%
  Applications of the change-of-rings spectral sequence to the 
  computation of Hochschild cohomology}
\author{%
  Mariano Su\'arez-Alvarez}
\address{%
  Facultad de Ciencias Exactas y Naturales\\
  Universidad de Buenos Aires\\
  Ciudad Universitaria, Pabell\'on I\\
  Buenos Aires (1428) Argentina.}
\email{mariano@dm.uba.ar}
\thanks{This work has been done as part of the research projects
\textsc{UBACyT}~X169, \textsc{PIP} \textsc{CoNICET}~5099 and
\textsc{PICS} \textsc{CNRS}}


\date{April 18, 2007}

\dedicatory{Dedicated to Henri Cartan and Samuel Eilenberg.}

\begin{abstract}
We consider the change-of-rings spectral sequence as it applies to
Hochschild cohomology, obtaining a description of the differentials on
the first page which relates it to the multiplicative stucture on
cohomology. Using this information, we are able to completely describe
the cohomology structure of monogenic algebras as well as some
information on the structure of the cohomology in more general
situations. 

We also show how to use the spectral sequence to reprove and
generalize results of M.~Auslander~\etal~about homological
epimorphisms. We derive from this a rather general version of the long
exact sequence due to D.~Happel for a one-point (co)-extension of a
finite dimensional algebra and show how it can be put to use in
concrete examples. \end{abstract}

\maketitle

\section*{Introduction}
\label{sect:intro}

The computation of the Hochschild cohomology~$\HH^\*(A)$ of an
algebra~$A$ is usually a difficult, laborious task. In most cases, the
well-known identification $\HH^\*(A)\cong\Ext_{A^e}^\*(A,A)$ and the
method developed by Cartan and Eilenberg in \emph{Homological
Algebra}~\cite{CartanEilenberg} for the computation of derived
functors are used: one finds an $A^e$-projective resolution $P^\*$ for
$A$ and then computes the cohomology of the
complex~$\hom_{A^e}(P^\*,A)$. Up to the flexibility of being able to
substitute the standard Hochschild resolution for a more convenient
one, this method is a direct implementation of the definition
of~$\HH^\*(A)$ given by Hochschild in~\cite{Hochschild}. 

As it is well known, the Hochschild cohomology~$\HH^\*(A)$ is endowed
naturally with both an associative algebra structure and a graded Lie
algebra structure. Recently, these structures have received
considerable attention: the r\^ole played by the former in the
representation theory of~$A$ is being studied by many researchers and
a theory modeled on Quillen's theory of support varieties for
groups~\cite{Quillen} is emerging, while the latter has been central in
recent developements related to the deformation theory of algebras.

Now, as soon as one attempts to make these structures explicit in
specific examples, computational difficulties arise very quickly. The
cup product can be computed using an arbitrary $A^e$-projective
resolution~$P^\*$ of~$A$ provided with a diagonal map $\Delta:P^\*\to
P^\*\ot_AP^\*$ but, while this is in general considerably more
convenient than dealing with the standard Hochschild resolution, this
quickly becomes impractical. The Lie structure, on the other hand, is
defined in terms of the Hochschild resolution and we do not have
available\dash---to the author's knowledge\dash---a way of computing
it in terms of an arbitrary resolution. Consequently, when a
projective $A^e$-resolution~$P^\*$ of~$A$ is being used to
determine~$\HH^\*(A)$, in order to compute the Lie structure one needs
comparison maps $P^\*\rightleftarrows B^\*(A)$ between $P^\*$ and the
Hochschild resolution~$B^\*(A)$. This is extremely messy;
in most cases, in fact, this approach is used successfully only
in low degrees.

This rather unsatisfactory situation should be compared to what
happens in the context of group cohomology. Indeed, while by no means
a trivial task, the effective computation of group cohomology as an
algebra is a rather well under\-stood process in which the various pieces
of structure present\dash---the Lyndon-Hochschild-Serre spectral
sequence, the action of the Steenrod algebra and, more generally, the
close relationship with algebraic topology, etc.\dash---serve as
powerful tools. At present, we do not have comparable tools
at hand when dealing with Hochschild cohomology.

\medskip

In this paper, we pick one of the general constructions of homological
algebra, the spectral sequence for a change of rings, and we try to
put it to use in the problem of computing Hochschild cohomology. It
turns out that the differentials on the initial term of this spectral
sequence can be described in terms of the Yoneda product. Our main
technical result, and the objective of sections~\pref{sect:operations}
and~\pref{sect:change-of-rings} below, is the following theorem:

\begin{Theorem*}
Let $k$ be a field and let $\phi:A\to B$ be a morphism of
$k$\nbd-algebras. If $M$ is a $B$\nbd-bimodule, there is a natural convergent
spectral sequence such that
  \[
  E_2^{p,q}\cong\Ext_{B^e}^p(\Tor_q^A(B,B),M)
  \Rightarrow
  H^\*(A,M).
  \]
Moreover, for each $q\geq1$ there exists
$\O_q(\phi)\in\Ext_{B^e}^2(\Tor_{q-1}^A(B,B),\Tor_q^A(B,B))$ such that
the differential on the term $E_2$,
  \[
  d_2^{p,q}:\Ext_{B^e}^p(\Tor_q^A(B,B),M)\to\Ext_{B^e}^{p+2}(\Tor_{q-1}^A(B,B),M),
  \]
is given by 
  \[
  d_2^{p,q}(\zeta)=(-1)^p\zeta\circ\O_q(\phi)
  \]
if $\zeta\in\Ext_{B^e}^p(\Tor_q^A(B,B),M)$ and $\circ$ is the Yoneda
product.
\end{Theorem*}

\medskip

We remark that this theorem provides invariants~$\O_q(\phi)$ for the
morphism~$\phi$. When $\phi$ is surjective, $\O_1(\phi)$~is
well-known:

\begin{Proposition*}
Assume, in the theorem, that $\phi:A\to B$ is surjective and put
$I=\ker\phi$. Then the class $\O_1(\phi)$ can be seen as an element of
$H^2(B, I/I^2)$ and then coincides with the characteristic class of the
singular extension of algebras
  \[\xymatrix{
  0 \ar[r]
    & I/I^2 \ar[r]
    & A/I^2 \ar[r]
    & B \ar[r]
    & 0
  }\]
\par\removelastskip
\end{Proposition*}

\medskip

The classes $\O_q(\phi)$ can be considered, then, as higher
degree analogues of the characteristic class. It would be interesting
to have tractable descriptions for them and their possible
interrelations. 

We remark that we have not been able to obtain information on the
differentials $d^p$ for $p>2$ of the spectral sequence appearing in
the theorem. Such information would be quite useful.

Using our theorem, in section~\pref{sect:monogenic}, we are able to
work out a computation of the Hochschild cohomology of algebras
$k[X]/(f)$ which are quotients of a polynomial algebra in one
variable, arriving at a presentation of the cohomology algebra
$\HH^\*(k[X]/(f))$ and of the Gerstenhaber Lie bracket on it; see
theorems~\pref{thm:monogenic} and~\pref{thm:monogenic:lie} below for
precise statements. Interestingly, this computation relies on explicit
work with resolutions only in very low degrees. The sort of arguments
used for this are of wider applicability: in
section~\pref{sect:variations} we present some variations and we
intend, in future work, to provide further related results.

Finally, in section~\pref{sect:nice} we consider the
change-of-rings spectral sequence corresponding to an homological
epimorphism $\phi:A\to B$, \ie, a morphism of algebras which induces a
full and faithful embedding $D^b(\lMod B)\to D^b(\lMod A)$ of bounded
derived categories. We are able to generalize results which relate the
Hochschild cohomologies of~$A$ and~$B$ due to
M.~Auslander, M.~I.~Platzeck and G.~Todorov~\cite{APT} and J. A. de la
Pe\~na and C.~Xi~\cite{delaPenhaXi}. In particular,
following~\cite{delaPenhaXi} we obtain long exact sequences
generalizing the one constructed by D.~Happel in~\cite{Happel} for 
one-point (co)extensions of finite dimensional algebras as well as
the generalizations of C.~Cibils~\cite{Cibils1}, E.~Green and
\oldO.~Solberg~\cite{GreenSolberg}, E.~Green, E.~N.~Marcos and
N.~Snashall~\cite{GreenMarcosSnashall}.

\paragraph*[Global conventions.] We fix a field $k$ throughout the
paper. Algebras will be always be $k$\nbd-algebras, $\hom$ and
$\ot$ will be taken over $k$ and, in general, all our linear
constructions will be $k$\nbd-linear. 

When working with algebras given as admissible quotients of path algebras,
our notations and nomenclature should be standard. As we prefer
left modules, we compose arrows from right to left. We consider paths in a
quiver as elements both of the path algebra, as usual, and of the quotient
algebras thereof, provided they represent a non-zero element in such a
quotient.

We refer to Stanley's book~\cite{Stanley}*{Chapter~3} for the little
we need about partially ordered sets.

\medskip

The author would like to thank Andrea Solotar and Estanislao
Herscovich heartily for their very careful reading of this manuscript.


\section{Stable operations on cohomology}
\label{sect:operations}

\paragraph\label{p:Yoneda} 
Let $A$ be a algebra. Recall that composition of extensions gives, 
for left $A$\nbd-modules $M$, $N$ and $P$, a Yoneda product
  \[
  \circ:\Ext_A^\*(N,P)\ot\Ext_A^\*(M,N)\to\Ext_A^\*(M,P)
  \]
which is homogeneous, natural and associative. Given a short exact
sequence $E:N'\mono N\epi N''$ of left $A$\nbd-modules and a
left $A$\nbd-module $M$, the connecting homomorphisms in the
long exact sequence for the $\partial$\nbd-functor
$\Ext_A^\*(M,-)$, namely
$\partial:\Ext_A^q(M,N'')\to\Ext_A^{q+1}(M,N')$, are given, for
$\sigma\in\Ext_A^q(M,N'')$, by the simple formula
  \[
  \partial\sigma=E\circ\sigma,
  \]
where we regard $E\in\Ext_A^1(N'',N')$. This is explained in detail in
\cite{MacLane}, chapter \textsc{III}.

\paragraph Given a pair of functors $F$, $T:\category{A}\to\category{B}$
between categories $\category{A}$ and $\category{B}$, we write $[F,T]$
the class of natural transformations $F\to T$. If $f\in[F,T]$ and
$x\in\category{A}$, we write $f_x:Fx\to Tx$ the component morphism of
$\category{B}$ corresponding to $x$.

\paragraph Fix now an algebra $A$. If $p$, $q\geq 0$ and $M$ and $N$
are left $A$\nbd-modules, set
  \[
  \Op_A^{p,q}(M,N)=[\,\Ext_A^p(N,-),\Ext_A^q(M,-)].
  \]
Clearly $\Op_A^{p,q}$ is an additive functor, contravariant on its first
variable, covariant on the second one. Yoneda's lemma gives us an
isomorphism of functors
  \[
  \Op_A^{0,q}(M,N)=[\,\hom_A(N,-),\Ext_A^q(M,-)]\cong\Ext_A^q(M,N).
  \]

\paragraph Given left $A$\nbd-modules $M$ and $N$ and
$d\geq0$, a \definiendum{cohomological operation $\O$ of degree $d$ from
$M$ to $N$} is a sequence $\O=(\O^p)_{p\geq0}$ of natural
transformations $\O^p\in\Op_A^{p,p+d}(M,N)$. We say that $\O$ is
\definiendum{stable} if, for each short exact sequence $P'\mono P\epi
P''$ of left $A$\nbd-modules, the following diagram commutes
for each $p\geq 0$:
  \[\xymatrix{
  \Ext_A^p(N,P'')\ar[r]^-\partial\ar[d]_-{\O^p_{P''}}
    & \Ext_A^{p+1}(N,P')\ar[d]^-{\O^{p+1}_{P'}} \\
  \Ext_A^{p+d}(M,P'')\ar[r]^-\partial
    & \Ext_A^{p+1+d}(M,P')
  }\]
We write $\sOp_A^d(M,N)$ the class of all stable cohomological operations
from~$M$ to~$N$. Again, it is clear that this is an additive functor, with
the same variances as $\Op_A^{\*,\*}(\place,\place)$.

\paragraph The Yoneda product allows us to define, for each $d\geq0$,
a natural morphism $Y:\Ext_A^d(M,N)\to\sOp_A^d(M,N)$ in the following way: if
$\zeta\in\Ext_A^d(M,N)$, $p\geq0$ and $P$ is a left
$A$\nbd-module, we let
$Y(\zeta)_P^p=(-)\circ\zeta:\Ext_A^p(N,P)\to\Ext_A^{p+d}(M,P)$ be
given by post-multiplication by $\zeta$. This is well defined because
$Y(\zeta)$ is a cohomological operation by definition, which is
stable because the Yoneda product is associative and the connecting
homomorphisms can be expressed as in~\pref{p:Yoneda}. The morphism~$Y$
is additive and natural both in~$M$ and~$N$.

It is a monomorphism; in fact, writing
$1_N\in\hom_A(N,N)=\Ext_A^0(N,N)$ the identity map of $N$, if
$\zeta\in\Ext_A^d(M,N)$, then $Y(\zeta)_N^0(1_N)=\zeta$, so we 
recover~$\zeta$ from~$Y(\zeta)$. On the other hand, if $\O\in\sOp_A^d(M,N)$
and $E_\O=\O^0_N(1_N)\in\Ext_A^d(M,N)$, put $\R=Y(E_\O)-\O$; this is a
stable cohomological operation of degree $d$. From the very definition,
$\R_N^0(1_N)=0$; if $f\in\hom_A(N,P)=\Ext_A^0(N,P)$ and, for each left
$A$\nbd-module $Q$, $f_Q^q:\Ext_A^q(Q,N)\to\Ext_A^q(Q,P)$ is the induced
morphism, naturality implies that
$\R_P^0(f)=\R_P^0f_N^0(1_N)=f_M^q\R_N^0(1_N)=0$. This means that in fact
$\R^0=0$. We will show in~\pref{p:R=0} that any stable cohomological
operation which vanishes---as $\R$ does---on $\Ext_A^0$ is identically
zero; this will allow us to conclude that $Y(E_\O)=\O$, proving the
following theorem:

\begin{Theorem}\label{thm:sOp=Ext}
There is an isomorphism of graded bifunctors
  \[
  Y : \Ext_A^\*(-,-)\xrightarrow{\cong}\sOp_A^\*(-,-).
  \]
\par\removelastskip
\end{Theorem}

\paragraph\label{p:R=0} Consider, then, a stable cohomological
operation $\O\in\sOp_A^d(M,N)$ and suppose that $p\geq0$ and
$\O^p=0$. Let $P$ be a left $A$\nbd-module and choose any
short exact sequence $P\mono I\epi P'$ in which $I$ is an injective
module. Stability of $\O$ entails the commutation of the following
diagram, in which the top row is exact:
  \[\xymatrix{
  \Ext_A^p(N,P')\ar[r]^-\partial\ar[d]_-{\O^p_{P'}}
    & \Ext_A^{p+1}(N,P)\ar[d]^-{\O^{p+1}_{P}}\ar[r]
    & \Ext_A^{p+1}(N,I)=0 \\
  \Ext_A^{p+d}(M,P')\ar[r]^-\partial
    & \Ext_A^{p+1+d}(M,P)\ar[r]
    & \Ext_A^{p+1+d}(M,I)=0
  }\]
The hypothesis that $\O^p=0$ implies that $\O_P^{p+1}=0$. The
arbitrariness of $P$ and an evident inductive argument show that
$\O=0$, as we needed.

\paragraph Our definition of stable cohomological operation restricts
us to consider operations of non-negative degree. The proof of
theorem~\pref{thm:sOp=Ext} shows that there is no loss in this, since
an operation of negative degree will certainly vanish on $\Ext_A^0$.
We note however that in general one has $\Op_A^{p,q}(M,N)\not=0$
for $p>q$, \cf~\cite{HiltonRees}.

\paragraph It is clear that under the isomorphism of the theorem, the
Yoneda product is identified with the composition of stable
operations.

\section{Change of rings}
\label{sect:change-of-rings}

\subsection{The spectral sequence}
\label{subsect:change}

\paragraph Consider a morphism of algebras $\phi:A\to B$ and left
$A$\nbd- and $B$\nbd-modules~$M$ and~$N$,
respectively. Let $X^\*\epi M$ be a projective resolution of~$M$ as an
$A$\nbd-module and $N\mono Y^\*$ an injective resolution of~$N$ 
as a $B$\nbd-module. We consider the complex 
$Z^\*=\hom_B(B\ot_AX^\*,Y^\*)$.

\paragraph The filtration $\I F^\*Z^\*$ on $Z^\*$ with
  \[
  \I F^pZ^q=\bigoplus_{\substack{r+s= q\\r\geq p}}
    \hom_B(B\ot_AX^r,Y^s)
  \]
determines a spectral sequence~$\I E$ converging to~$H(Z^\*)$. The
differential on~$\I E_0$ is induced by the one on~$Y^\*$ and, since
there is an isomorphism of functors $\hom_B(B\ot_A\place,\place)\cong
\hom_A(\place,\place)$ and each $X^p$, when $p\geq0$, is a projective
$A$\nbd-module, we have that $\I E_1^{p,q}=0$ if $q>0$ and $\I
E_1^{\*,0}\cong\hom_A(X^\*,N)$. The differential on $\I E_1$
corresponds in a natural way to the differential on $X^\*$, so that
$\I E_2^{\*,0}\cong\Ext_A^\*(M,N)$. The spectral sequence degenerates
at the term $\I E_2$ and convergence implies that
$H(Z^\*)\cong\Ext_A^\*(M,N)$.

\paragraph Considering now the filtration $\II F^\*Z^\*$ given by
  \[
  \II F^qZ^p=\bigoplus_{\substack{r+s=p\\s\geq q}}
    \hom_B(B\ot_AX^r,Y^s)
  \]
we obtain another cohomologically graded spectral sequence contained
in the first quadrant converging to $\Ext_A^\*(M,N)$. We have $\II
E_0^{p,q}\cong\hom_B(B\ot_A X^q,Y^p)$ with differential induced by the
one on $X^\*$; as $\hom_B(-,Y^p)$ is an exact functor for each $p\geq
0$, we find that $\II E_1^{p,q}\cong\hom_B(\Tor^A_q(B,M),Y^p)$. The
differential on the term $\II E_1$ is induced by the differential on
$Y^\*$, so that the next term has $\II
E_2^{p,q}\cong\Ext_B^p(\Tor^A_q(B,M),N)$. 

We record these facts in the following proposition.

\begin{Proposition}\label{prop:change}
Let $\phi:A\to B$ be a morphism of algebras. Let $M$ be a left
$A$\nbd-module and $N$ a left $B$\nbd-module and
consider $N$ as a left $A$\nbd-module by pull-back along
$\phi$. There is a convergent spectral sequence
  \[
  E_2^{p,q}\cong\Ext_B^p(\Tor^A_q(B,M),N)\Rightarrow\Ext_A^\*(M,N).
  \]
It is natural with respect to both $M$ and $N$.~\qed
\end{Proposition}

\paragraph This spectral sequence was first constructed
in~\cite{CartanEilenberg}*{\textsc{XVI}.5, case 3}.

\paragraph Standard properties of resolutions and the comparison
theorem for spectral sequences immediately imply that the spectral
sequence of~\pref{prop:change} does not depend on the particular
resolutions $X^\*$ and $Y^\*$ used to construct them.

\paragraph\label{p:edge} The spectral sequence constructed
in~\pref{subsect:change} is in the first quadrant, so it comes
equipped with an edge morphism from the limit to the ``base,''
  \[
  e:\Ext_A^\*(M,N) \to E_2^{0,\*}\cong\hom_B(\Tor_\*^A(B,M),N)
  \]
and a morphism from the ``fiber'' to the limit,
  \begin{equation}\label{eq:fiber}
  e':E_2^{\*,0}\cong\Ext_B^\*(B\ot_AM,N)\to\Ext_A^\*(M,N).
  \end{equation}
They can be computed as follows.

Let $X^\*\epi M$ be a projective resolution of $M$ as a left
$A$\nbd-module and take a class $\alpha\in\Ext_A^p(M,N)$. To compute
$e(\alpha)$, we pick a representative $a\in\hom_A(X^p,N)$ for~$\alpha$ and
let $\bar a\in\hom_B(B\ot_A X^p,N)$ be the morphism corresponding to it by
the natural identification, so that $\bar a(b\ot x)=ba(x)$. If
$\tau\in\Tor_p^A(B,M)$ is represented by $t\in B\ot_AX^p$ in the complex
$B\ot_AX^\*$, then
  \[
  e(\alpha)(\tau) = \bar a(t).
  \]
Consider now additionally a projective resolution $X'^\*\epi B\ot_AM$
of $B\ot_AM$ as a left $B$\nbd-module. There exists a morphism of
complexes of left $B$\nbd-modules $f:B\ot_AX^\*\to X'^\*$ over
the identity map of $B\ot_AM$, well determined up to homotopy;
indeed, the graded components of $B\ot_AX^\*$ are
$B$\nbd-projective. The morphism $e'$ is then induced on homology
by the composition
  \[
  \hom_B(X'^\*,N)
  \xrightarrow{f^*}
  \hom_B(B\ot_A X^\*,N)
  \cong
  \hom_A(X^\*,N).
  \]

\paragraph\label{p:e'-with-exts} Suppose that $\phi:A\to B$ is a
surjection of algebras. Then, of course, there is a canonical
isomorphism $B\ot_AM\cong M$ if $M\in\lMod B$ and the edge morphism
from the fiber to the limit in the spectral sequence is, when both $M$
and $N$ are left $B$\nbd-modules,
  \[
  e':\Ext_B^\*(M,N)\to\Ext_A^\*(M,N).
  \]
In this situation, we can describe $e'$ in terms of iterated
extensions of modules: if~$\alpha\in\Ext_B^p(M,N)$ is the class of
the $p$\nbd-extension of left $B$\nbd-modules
  \begin{equation}\label{eq:p-ext}
  \xymatrix{
  0  \ar[r]
    & N \ar[r]
    & U_p \ar[r]
    & \cdots \ar[r]
    & U_1 \ar[r]
    & M \ar[r]
    & 0
  }\end{equation}
then $e'(\alpha)\in\Ext_A^p(M,N)$ is simply the
$p$\nbd-extension~\eqref{eq:p-ext} now considered \emph{in}~$\lMod A$.
This follows easily from the recipe given in~\pref{p:edge} for
computing~$e'$ and the details of the proof that the
$\Ext$\nbd-functors can be computed from iterated extensions as presented
in~\cite{HiltonStammbach}*{Theorem~\textsc{IV}.9.1}. We leave this to
the reader and just recall briefly how one goes from iterated
extensions to cocycles, as we will need this below.

Let $M\in\lMod A$ and fix a projective resolution
$X_\*\to M$ of $M$ in $\lMod A$. Consider a $p$\nbd-extension $\xi$ in that
category ending in $M$, as in the bottom row of the following diagram
  \[\xymatrix@C-5pt{
  {}
    & \cdots\ar[r]
    & X_p \ar[r] \ar@{..>}[d]^-{h_p}
    & X_{p-1} \ar[r] \ar@{..>}[d]^-{h_{p-1}}
    & \cdots \ar[r] 
    & X_1 \ar[r] \ar@{..>}[d]^-{h_0}
    & X_0 \ar[r] \ar@{..>}[d]^-{h_0}
    & M \ar[r] \ar[d]^-{\id_M}
    & 0
    \\
  \xi:
    & 0  \ar[r]
    & N \ar[r]
    & U_p \ar[r]
    & \cdots \ar[r]
    & U_2 \ar[r]
    & U_1 \ar[r]
    & M \ar[r]
    & 0
  }\]
Since the top row is made of projective modules and the bottom one is
exact, the identity map $\id_M$ on the right extends to a map of
complexes $h_\*$. Clearly $h_p:X_p\to N$ is actually a
$p$\nbd-cocycle in $\hom_A(X_\*,N)$; $h_\*$ is defined only up
to homotopy and this means that $h_p$ is defined up to a coboundary in
$\hom_A(X_\*,N)$. The class of $h_p$ in $H(\hom_A(X_\*,N))$ depends
thus only on $\xi$ and it is its characteristic class in
$\Ext_A(M,N)$.

\paragraph This description~\pref{p:e'-with-exts} of the
morphism~\eqref{eq:fiber} in terms of extensions has the immediate
consequence that $e'$ is multiplicative on $B$\nbd-modules: if
$M$, $N$, $P\in\lMod B$ and $\alpha\in\Ext_B^p(M,N)$,
$\beta\in\Ext_B^q(N,P)$, then
  \[
  e'(\beta\circ\alpha) = e'(\beta)\circ e'(\alpha),
  \]
if $\circ$ denotes Yoneda composition of iterated extensions.

\subsection{The differentials}
\label{subsect:differentials}

\paragraph In the following lemma, we consider double complexes $X^{\*,\*}$
whose horizontal and vertical differentials $\delta'$ and $\delta''$
anti-commute and cohomologically graded spectral sequences~$E^{\*,\*}_\*$
for which the first upper index corresponds to the filtration degree.

\begin{Lemma*}\label{lemma:tech}
Let 
  \[\xymatrix{
  0 \ar[r]
    & \pre1X{^{\*,\*}} \ar[r]^-{j_0}
    & \pre2X{^{\*,\*}} \ar[r]^-{k_0}
    & \pre3X{^{\*,\*}} \ar[r]
    & 0
  }\]
be a short exact sequence of double complexes such that the induced
sequence 
  \[\xymatrix{
  0 \ar[r]
    & \pre1E{_1^{\*,\*}} \ar[r]^-{j_1}
    & \pre2E{_1^{\*,\*}} \ar[r]^-{k_1}
    & \pre3E{_1^{\*,\*}} \ar[r]
    & 0
  }\]
is also exact. If $\partial:\pre3E{_2^{p,q}}\to\pre1E{_2^{p,q+1}}$ is
the connecting homomorphism corresponding to the differentials in this
last sequence, the square
  \[\xymatrix{
  {}\pre3E{_2^{p,q}}\ar[r]^-\partial\ar[d]_-{d_2}
    & {}\pre1E{_2^{p,q+1}}\ar[d]^-{d_2} \\
  {}\pre3E{_2^{p-1,q+2}}\ar[r]^-\partial
    & {}\pre1E{_2^{p-1,q+3}}
  }\]
anti-commutes.
\end{Lemma*}

\medskip

\begin{figure}
  \vskip-30pt 
  \[\begin{xy}
  \xymatrix"a"@!=2pt{
    \\
    d_2\alpha \\
    b\ar[r]\ar[u]
      & \cdot \\
    *+{t}\ar@{..}[u]\ar@{.}[r]
      & a\ar[u]\ar[r]
      & 0
    }
  \POS*{}
  \POS-(25,15)
  \xymatrix"b"@!=2pt{
    \\
    \cdot\ar[r]\ar@{..}["a"]
      & \cdot\\
    s\ar[u]\ar@{.}[r]\ar@{..}["a"]="q"
      & \cdot\ar@{..}["a"]\ar[u]\\
    {} 
      & c\ar[u]\ar[r]\ar@{-->}["a"]^(.50){}="c"
      & 0
    }
  \POS*{}
  \POS-(25,15)
  \xymatrix"c"@!=2pt{
    d_2\partial\alpha\\
    r\ar[u]\ar[r]\ar@{..}["b"]
      & \cdot\ar["b"]\\
    {} 
      & x\ar[u]\ar[r]\ar@{-->}["b"]^(.50){}="x"
      & 0
    }
  {
  \SelectTips{cm}{10}
  \POS "b3,1" \ar@/l/@{-*{\UseTips\dir{>}}} "x"
  \POS "a4,1" \ar@/l/@{-*{\UseTips\dir{>}}} "c"
  }
  \POS "q"*{}
  \POS(35,-20)
  \xymatrix@C-5pt@R-5pt{
    & \\
    *{}=0\ar[u]^-{\delta''}^>q\ar[r]_-{\delta'}_>p&
    }
  \end{xy}
  \]
\caption{The chase in the proof of~\pref{lemma:tech}.}\label{fig:chase}
\end{figure}
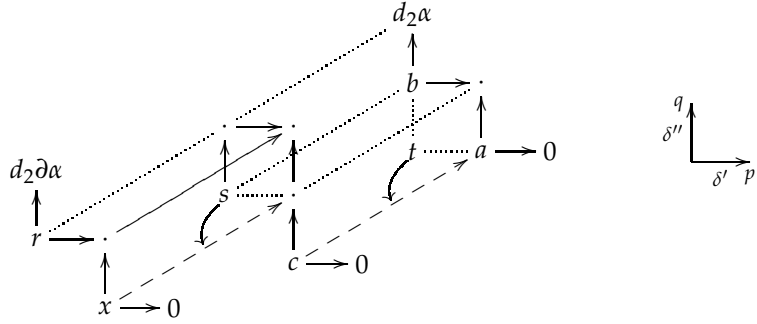

The diagram in figure~\ref{fig:chase} may be of help in following the proof
of the lemma. The dotted lines are used to show the relative positions of
the elements appearing in the diagram, solid and broken arrows represent
maps on $E_0$ and $E_1$, respectively, and the curved lines show that, for
example, $j_0x=\delta''c+\delta's$; finally, the planes are, in order of
increasing depth, $\pre1E{_0}$, $\pre2E{_0}$ and~$\pre3E{_0}$.

\begin{proof}
Let $\alpha\in\pre3E{_2^{p,q}}$ and $a\in\pre3E{_0^{p,q}}$ be such
that $a\in\alpha$; then $\delta'a=0$ and there exists
$b\in\pre3E{_0^{p-1,q+1}}$ such that $\delta'b=\delta''a$. We see that
$\delta''b\in d_2\alpha$. As $k_1$ is an epimorphism, there exists
$c\in\pre2E{_0^{p,q}}$ and $t\in\pre3E{_0^{p-1,q}}$ such that
$\delta'c=0$ and $k_0c=a+\delta't$. We have
  \[
  k_0\delta''c=\delta''k_0c=\delta''a+\delta''\delta't=\delta'(b-\delta''t),
  \]
and since $\ker k_1=\im j_1$, there are $x\in\pre1E{_0^{p,q+1}}$ and
$s\in\pre2E{_0^{p-1,q+1}}$ with $\delta'x=0$ and
$j_0x=\delta''c+\delta's$. Note that $x\in\partial\alpha$.

Observe now that $j_0\delta''x= \delta''j_0x= \delta''\delta's=
-\delta'\delta''s$; as $j_1$ is a monomorphism, there is
$r\in\pre1E{_0^{p-1,q+2}}$ such that $\delta'r=\delta''x$. Then
$\delta''r\in d_2\partial\alpha$.

Let now $z=j_0r+\delta''s\in\pre2E{_0^{p-1,q+2}}$. We have that
  \[
  \delta'k_0s=k_0\delta's=k_0j_0x-k_0\delta''c=-\delta''k_0c
    =-\delta''a-\delta''\delta't=-\delta''a+\delta'\delta''t,
  \]
so $\delta'(\delta''t-k_0s)=\delta''a$. This implies that 
  \[
  k_0z = k_0\delta''s 
       = \delta''k_0s 
       = \delta''(k_0s-\delta''t)\in-d_2\alpha
  \]
and because $j_0\delta''r=\delta''j_0r=\delta''z$,
$\delta''r\in-\partial d_2\alpha$. We thus see that
$d_2\partial\alpha=-\partial d_2\alpha$, as we were required to show.
\end{proof}

\paragraph\label{p:mess} We will apply the lemma in the context of the
spectral sequence from~\pref{prop:change}.  Consider a morphism of
algebras $\phi:A\to B$, a left $A$\nbd-module $M$ and a short
exact sequence of left $B$\nbd-modules 
  \begin{equation}\label{eq:suc:N}
  \xymatrix{
  0 \ar[r]
    & \pre1N{} \ar[r]^-j
    & \pre2N{} \ar[r]^-k
    & \pre3N{} \ar[r]
    & 0
  }
  \end{equation}
Let $X^\*\epi M$ be a projective resolution of $M$, and let us choose
injective resolutions ${}_iN\mono{}_iY^\*$, for $1\leq i\leq 3$, and
morphisms $j^\*$ and $k^\*$ among these such that in the diagram
  \[\xymatrix{
  0\ar[r]
    & {}\pre1N{}\ar[r]^-j\ar[d]
    & {}\pre2N{}\ar[r]^k\ar[d]
    & {}\pre3N{}\ar[r]\ar[d]
    & 0\\
  0\ar[r]
    & {}\pre1Y{^\*}\ar[r]^-{j^\*}
    & {}\pre2Y{^\*}\ar[r]^-{k^\*}
    & {}\pre3Y{^\*}\ar[r]
    & 0
  }\]
the rows are exact and each square commutes. If we set
  \[
  \pre iZ{^{\*,\*}}=\hom_B(B\ot_A X^\*,\pre iY{^\*})
  \]
for each $i$ with $1\leq i\leq 3$, we have, as one can easily show, an
exact sequence of double complexes 
  \[\xymatrix{
  0 \ar[r]
    & \pre1Z{^{\*,\*}} \ar[r]
    & \pre2Z{^{\*,\*}} \ar[r]
    & \pre3Z{^{\*,\*}} \ar[r]
    & 0
  }\]
in which the morphisms are induced by $j^\*$ and $k^\*$. We consider
spectral sequences as in proposition~\pref{prop:change}.

For each bidegree, the exact sequence
  \[\xymatrix{
  0 \ar[r]
    & \pre1E{_0^{p,q}} \ar[r]
    & \pre2E{_0^{p,q}} \ar[r]
    & \pre3E{_0^{p,q}} \ar[r]
    & 0
  }\] 
is obtained by applying the functors $\hom_B(B\ot_A X^q,-)$ to
  \[\xymatrix{
  0 \ar[r]
    & \pre1Y{^p} \ar[r]
    & \pre2Y{^p} \ar[r]
    & \pre3Y{^p} \ar[r]
    & 0
  }\]
Since this last sequence splits, taking homology we get an exact
sequence 
  \[\xymatrix{
  0 \ar[r]
    & \pre1E{_1^{p,q}} \ar[r]
    & \pre2E{_1^{p,q}} \ar[r]
    & \pre3E{_1^{p,q}} \ar[r]
    & 0
  }\]
We are then in the situation of the lemma and we see that the
following diagram is anti-commutative:
\begingroup
\small
  \[
  \mskip-30mu
  \xymatrix@C-5pt{
  {}\pre3E{_2^{p,q}}\cong\Ext_B^p(\Tor_q^A(B,M),\pre3N{})
      \ar[r]^-\partial\ar[d]_-{d_2^{p,q}} 
    & {}\pre1E{_2^{p+1,q}}\cong\Ext_B^{p+1}(\Tor_q^A(B,M),\pre1N{})
      \ar[d]^-{d_2^{p+1,q}} \\
  {}\pre3E{_2^{p+2,q-1}}\cong\Ext_B^{p+2}(\Tor_{q-1}^A(B,M),\pre3N{})
      \ar[r]^-\partial
    &
    {}\pre1E{_2^{p+3,q-1}}\cong\Ext_B^{p+3}(\Tor_{q-1}^A(B,M),\pre1N{})
  }\]
\endgroup
where $\partial:\pre3E{_2^{p,q}}\to\pre1E{_2^{p+1,q}}$ is the
connecting homomorphism for the long exact sequence corresponding to
the $\partial$\nbd-functor $\Ext_B^\*(\Tor_q^A(B,M),\place)$
and the short exact sequence~\eqref{eq:suc:N}.

We are thus led to the following theorem:

\begin{Theorem}\label{thm:diffs}
Let $\phi:A\to B$ be a morphism of algebras and let $M$ be a left
$A$\nbd-module. For each $q\geq0$ there is a class
$\O_q(M)\in\Ext_B^2(\Tor_{q-1}^A(B,M),\Tor_q^A(B,M))$ such that, for
each left $B$\nbd-module $N$, the $E_2$ term of the spectral
sequence of proposition~\pref{prop:change} has differentials 
  \(
  d_2^{p,q}:\Ext_B^p(\Tor_q^A(B,M),N)\to\Ext_B^{p+2}(\Tor_{q-1}^A(B,M),N)
  \)
given by $d_2^{p,q}(\zeta)=(-1)^p\;\zeta\circ\O_q(M)$.
\end{Theorem}

\begin{proof}
The observations of~\pref{p:mess}, together with the naturality of the
spectral sequences involved, imply that if we define
$\O_q(M)_N^p=(-1)^pd_2^{p,q}$ for each left $B$\nbd-module~$N$, we obtain
an operation $\O_q(M)\in\sOp_B^2(\Tor_{q-1}^A(B,M),\Tor_{q}^A(B,M))$. The
theorem follows now from the description given in~\pref{thm:sOp=Ext} of
this set.
\end{proof}

\subsection{Hochschild cohomology}
\label{subsect:hochschild}

\paragraph Recall that the \definiendum{Hochschild
cohomology} of a $k$\nbd-algebra $A$ is the functor 
of $A$\nbd-bimodules
  \(
  H^\*(A,\place) = \Ext_{A^e}^\*(A,\place)
  \); here, as usual, $A^e= A\ot A^\op$ is the
so-called \definiendum{enveloping algebra} of $A$, which is such that
there is an isomorphism of categories $\lMod{A^e}\cong\bMod AA$. 

We write $\HH^\*(A)= H^\*(A,A)$ the value of Hochschild
cohomology on the regular $A$\nbd-bimodule $A\in\bMod AA$.

\paragraph In this context, theorem~\pref{thm:sOp=Ext} amounts to the
following: given $d\geq0$, there is a bijection between the class of
stable operations $H^\*(A,\place)\to H^{\*+d}(A,\place)$ and
$\HH^d(A)$ and these bijections can be collected into an algebra
isomorphism between the graded ring $\sOp_{A^e}^\*(A,A)$ of stable
operations on~$H^\*(A,\place)$ and the Yoneda algebra~$\HH^\*(A)$.

We remark that the Yoneda product on~$\HH^\*(A)$ coincides with the
classical cup product on Hochschild cohomology described
in~\cite{Gerstenhaber} and that the action of~$\HH^\*(A)$
on~$H^\*(A,\place)$ that corresponds under the isomorphism 
$\HH^\*(A)\cong\sOp_{A^e}^\*(A,A)$
to the action of~$\sOp_{A^e}^\*(A,A)$ 
on~$H^\*(A,\place)$ can itself be computed using cup products. 

\paragraph We want to make explicit the specialization
of~\pref{prop:change} and~\pref{thm:diffs} to Hochschild cohomology:

\begin{Theorem*}\label{thm:E:hoch}
Let $\phi:A\to B$ be a morphism of $k$\nbd-algebras. If
$M\in\bMod BB$, there is a convergent spectral sequence such that
  \[
  E_2^{p,q}\cong\Ext_{B^e}^p(\Tor_q^A(B,B),M)
  \Rightarrow
  H^\*(A,M)
  \]
and this spectral sequence is functorial on $M$. Moreover, for each
$q\geq0$ there exists
$\O_q(\phi)\in\Ext_{B^e}^2(\Tor_{q-1}^A(B,B),\Tor_q^A(B,B))$ such that
the differential of the term $E_2$,
  \[
  d_2^{p,q}:\Ext_{B^e}^p(\Tor_q^A(B,B),M)\to\Ext_{B^e}^{p+2}(\Tor_{q-1}^A(B,B),M),
  \]
is given by 
  \[
  d_2^{p,q}(\zeta)=(-1)^p\zeta\circ\O_q(\phi)
  \]
if $\zeta\in\Ext_{B^e}^p(\Tor_q^A(B,B),M)$.
\end{Theorem*}

\begin{proof}
The morphism $\phi$ induces in an obvious way a morphism
$\phi^e:A^e\to B^e$, to which one can apply~\pref{prop:change}. This
provides a spectral sequence with differentials of the form described
in~\pref{thm:diffs}. The limit of the spectral sequence is
$H^\*(A,M)$. On the other hand,
$\Tor_\*^{A^e}(B^e,A)\cong\Tor_\*^A(B,B)$ canonically,
\cf~\cite{CartanEilenberg}*{Corol.~\textsc{IX}.4.4} and this shows
that the $E_2$ is of the form stated.
\end{proof}

\paragraph This spectral sequence has good multiplicative properties. A
particularly useful one is the following:

\begin{Proposition*}\label{p:edge:cup}
Let $\phi:A\to B$ be a morphism of $k$\nbd-algebras such that
$B\ot_AB\cong B$. The edge morphism from the fiber to the limit
in the spectral sequence described in~\pref{thm:E:hoch} when $M$ is the
regular $B$\nbd-bimodule $B$ is then a map
  \[
  e':\HH^\*(B)\to H^\*(A,B).
  \]
This map is a morphism of algebras, when both its domain and its
codomain are endowed with the cup product.
\end{Proposition*}

\begin{proof}
This follows immediately from the description of~$e'$ given
in~\pref{p:edge} if one uses the standard Hochschild resolution in
order to compute $\Ext$\nbd-groups.
\end{proof}

\paragraph When the morphism $\phi$ is surjective, the class
$\O_1(\phi)$ is connected with a well-known construction:

\begin{Proposition*}\label{prop:O1}
Assume $\phi:A\to B$ is a surjection of algebras with kernel
$I=\ker\phi$. Then $I/I^2$ is a $B$\nbd-bimodule in a
natural way and there are isomorphisms $\Tor_0^A(B,B)\cong B$ and
$\Tor_1^A(B,B)\cong I/I^2$. Furthermore, the class
  \[
  \O_1(\phi)\in\Ext_{B^e}^2(\Tor_0^A(B,B),\Tor_1^A(B,B))= H^2(B,I/I^2)
  \]
constructed in~\pref{thm:E:hoch} coincides with the class of the
singular extension of algebras
  \begin{equation}\label{prop:O1:ext}
  \xymatrix{
  0 \ar[r]
    & I/I^2 \ar[r]
    & A/I^2 \ar[r]
    & B \ar[r]
    & 0
  }\end{equation}
\par\removelastskip
\end{Proposition*}

\begin{proof}
The existence of the claimed isomorphisms is left as an easy
exercise. The rest of the proposition can be proved by noting that
the exact sequence for terms of lower degree in the spectral sequence 
  \(
  E_2^{p,q}\cong\Ext_{B^e}^p(\Tor_q^A(B,B),\place)
  \Rightarrow
  H^\*(A,\place)
  \)
constructed in the theorem is the analog for Hochschild cohomology of
the $5$\nbd-term exact sequence constructed for group
cohomology in~\cite{HiltonStammbach}*{Section~\textsc{VI}.8} and then
adapting the arguments given in \cite{HiltonStammbach}*{Section~\textsc{VI}.10}.
\end{proof}

\paragraph We can construct explicitly a $2$-extension of $B$-bimodules
representing the class~$\O_1(\phi)$ appearing in~\pref{prop:O1}. The
construction is surely well-known but it does not appear in the standard
references.

Recall from~\cite{Karoubi} that if $\Lambda$ is an algebra, the
\definiendum{$\Lambda$-bimodule of non-commutative differential forms
$\Omega(A)$} is the kernel of the multiplication map
$\mu:\Lambda\ot\Lambda\to\Lambda$. In particular, there is an exact
sequence of $\Lambda$-bimodules
  \begin{equation}\label{eq:omega}
  \xymatrix{
  0 \ar[r]
    & \Omega(\Lambda) \ar[r]
    & \Lambda\ot\Lambda \ar[r]^-\mu
    & \Lambda \ar[r]
    & 0
  }\end{equation}
There is a map $d:\Lambda\to\Omega(\Lambda)$ given by
$d(\lambda)=\lambda\ot 1-1\ot\lambda$. One easily checks that this is a
derivation and, moreover, it turns out that $d$ is the universal derivation
of $\Lambda$ into $\Lambda$-bimodules, in the sense that composition with
$d$ induces an isomorphism $\hom_{\Lambda^e}(\Omega(\Lambda),\place) \cong
\Der(\Lambda,\place)$ of functors defined on $\Bimod\Lambda$.

Let us put ourselves back in the situation of~\pref{prop:O1}. Applying the
functor ${B^e\ot_{A^e}(\place)}$ to the short exact sequence~\eqref{eq:omega}
corresponding to $\Lambda=A$ and using the canonical isomorphism
$\Tor_1^{A^e}(B^e,A)\cong I/I^2$ noted above, we obtain an exact sequence
  \begin{equation}\label{eq:2ext}
  \xymatrix{
  0 \ar[r]
    & I/I^2 \ar[r]
    & B\ot_A\Omega(A)\ot_AB \ar[r]
    & B\ot B \ar[r]
    & B \ar[r]
    & 0
  }\end{equation}
where the second map is induced by $x\in I\mapsto 1\ot d(x)\ot 1\in
B\ot_A\Omega(A)\ot_AB$ and the following two are induced from the corresponding
maps in~\eqref{eq:omega}.

\begin{Proposition*}
The class~$\O_1(\phi)\in H^2(B,I/I^2)$ constructed in~\pref{prop:O1} is represented by the
$2$-extension of $B$-bimodules~\eqref{eq:2ext}.
\end{Proposition*}

\begin{proof}
Let $\sigma:B\to A/I^2$ be any $k$-linear section of the algebra map
$\bar\phi:A/I^2\to B$ induced by $\phi$. Define now
  \[
  c:b\ot b'\in B\ot B\mapsto \sigma(bb')-\sigma(b)\sigma(b')\in A/I^2.
  \]
Clearly $\bar\phi\circ c=0$, so there exists a unique $k$-linear map
$\alpha:B\ot B\to I/I^2$ such that $c$ is the composition of $\alpha$ and
the inclusion $I/I^2\hookrightarrow A/I^2$. A computation shows that
$\alpha$ is a $2$-cocycle in the standard Hochschild complex which computes
$H^\*(B,I/I^2)$ and, indeed, its class in $H^2(B,I/I^2)$ is the class
corresponding to the extension~\eqref{prop:O1:ext}. 

Let us pick any $k$-linear section $\iota:A/I^2\to A$ for the
projection $A\to A/I^2$ and define maps
  \[
  \hat\alpha:b\ot b'\ot b''\ot b'''\in B\ot B\ot B\ot B
        \mapsto b\alpha(b'\ot b'')b'''\in I/I^2
  \]
and
  \[
  \hat\sigma:b\ot b'\ot b''\in B\ot B\ot B      
        \mapsto b\ot d(\iota(\sigma(b')))\ot b''
        \in B\ot_A\Omega(A)\ot_AB;
  \]
that this is well defined follows from a simple computation. Then we have a
commutative diagram
  \[\xymatrix{
  \cdots \ar[r]
    & B^{\ot4} \ar[r] \ar[d]^-{\hat\alpha}
    & B^{\ot3} \ar[r] \ar[d]^-{\hat\sigma}
    & B\ot B \ar[r] \ar@{=}[d]
    & B \ar[r] \ar@{=}[d]
    & 0
    \\
  0 \ar[r]
    & I/I^2 \ar[r]
    & B\ot_A\Omega(A)\ot_AB \ar[r]
    & B\ot B \ar[r]
    & B \ar[r]
    & 0 
  }\]
in which the top row is the standard Hochschild resolution of~$B$.
Recalling the way in which one shows that Yoneda functors are computable
using projective resolutions, we see at once that this means that the
$2$-extension~\eqref{eq:2ext} represents the $2$-cocycle $\alpha$ and,
hence, the class $\O_1(\phi)$.
\end{proof}
 
\paragraph The classes $\O_q(\phi)\in
\Ext_{B^e}^2(\Tor_{q-1}^A(B,B),\Tor_q^A(B,B))$ constructed above for a
surjective algebra map $\phi:A\to B$ can be seen as higher order
invariants for the extension associated to~$\phi$,
extending\dash---according to~\pref{prop:O1}\dash---the classical invariant
$\O_1(\phi)$ constructed by Hochschild in~\cite{Hochschild}.

\paragraph Finally, we make explicit the $5$-term exact sequence to which
we made reference in the proof of~\pref{prop:O1} and describe the edge
morphism which appears in it. It is obtained as usual from the spectral
sequence constructed in~\pref{thm:E:hoch} by looking at terms of low
degree.

\begin{Proposition*}
Let $\phi:A\to B$ be a surjection of algebras and let $M$ be a
$B$-bimodule. There is an exact sequence
  \begin{multline*}
  \xymatrix@C-5pt{
  0 \ar[r]
    & H^1(B,M) \ar[r]
    & H^1(A,M) \ar[r]^-e
    & \hom_{B^e}(I/I^2, M) \ar[r] 
    &
    } \\
  \xymatrix@C-5pt{
  {} \ar[r]
    & H^2(B,M) \ar[r]
    & H^2(A,M)
  }
  \end{multline*}
The edge morphism $e$ is induced by the map $\Der(A,M) \to
\hom_{B^e}(I/I^2, M)$ which sends $f\in\Der(A,M)$ to the composition
  \[\xymatrix@C+10pt{
  I/I^2 \ar[r]
    & B\ot_A\Omega(A)\ot_AB \ar[r]^{1\ot f\ot 1}
    & B\ot_AM\ot_AB \ar[r]^-{\cong}
    & M
  }\]
The first map here is the one appearing in~\eqref{eq:2ext}.~\qed
\end{Proposition*}

\section{The cohomology of monogenic algebras}
\label{sect:monogenic}

\paragraph We fix a monic polynomial $f=\sum_{i=0}^N\alpha_i X^i\in
k[X]$ of degree $N$ and consider the algebra $A=k[X]/(f)$; let $x$ be
the class of $X$ in $A$. We want to describe its Hochschild cohomology
with values in an $A$\nobreakdash-bimodule $M$. To do so, we consider
the canonical projection $\phi:k[X]\to A$ and study the spectral
sequence $E$ attached to it in~\pref{thm:E:hoch}.

\paragraph The initial term $E_2$, which has
$E_2^{p,q}=\Ext_{A^e}^p(\Tor_q^{k[X]}(A,A),M)$, is easy to compute.
Indeed, the obvious short exact sequence
  \begin{equation}\label{eq:CCA}
  \xymatrix@1{
  0 \ar[r]
    & k[X]\ar[r]^-f
    & k[X]\ar[r]
    & A \ar[r]
    & 0
  }\end{equation}
provides a projective resolution of $A$ as a $k[X]$\nobreakdash-module
either on the left or on the right; using it to compute
$\Tor_\*^{k[X]}(A,A)$, we immediately see that
  \[
  \Tor_p^{k[X]}(A,A) \cong
    \begin{cases}
      A, & \text{if $p=0$ or $p=1$;} \\
      0, & \text{otherwise.}
    \end{cases}
  \]
This implies that there are isomorphisms
  \[
  E_2^{p,q} \cong
    \begin{cases}
      H^p(A,M), & \text{if $q=0$ or $q=1$;} \\
      0, & \text{otherwise.}
    \end{cases}
  \]

\paragraph We look now at the limit of $E$, which is $H^\*(k[X],M)$.
There is a projective resolution
  \[\xymatrix{
  0 \ar[r]
    & k[X]\ot k[X] \ar[r]^-{d_1}
    & k[X]\ot k[X] \ar@{->>}[r]^-{d_0}
    & k[X]
  }\]
of $k[X]$ as a $k[X]$\nobreakdash-bimodule, with $d_1(1\ot
1)=X\ot1-1\ot X$ y $d_0(1\ot1)=1$. A~trivial computation using it
shows that
  \[
  H^p(k[X],M) \cong
    \begin{cases}
      M^x, & \text{if $p=0$;} \\
      M_x, & \text{if $p=1$;} \\
      0, & \text{if $p\geq2$.} 
    \end{cases}
  \]
Here we are writing
  \(
  M^x \denotes \{m\in M:xm=mx\} 
  \)
and
  \(
  M_x \denotes \dfrac M{\{xm-mx:m\in M\}}
  \).

\paragraph We conclude that the $E_2$ term looks like this:
  \[
  \mskip-30mu
  \xymatrix@R-20pt@C-15pt{
  & \\
  & 0 & 0 & 0 & 0 & \cdots \\
  & {}H^0(A,M)\ar[rrd] & {}H^1(A,M)\ar[rrd] & {}H^2(A,M)\ar[rrd] 
   	& {}H^3(A,M) & \cdots \\
  & {}H^0(A,M) & {}H^1(A,M) & {}H^2(A,M) 
        & {}H^3(A,M) & \cdots \\
  *{}**{<0,0>}\ar[uuuu]^>q\ar[rrrrrrr]_>p&&&&&&&&\\
  & M^x & M_x & 0 & 0 & \cdots
  }\]
We have written the limit of this sequence, $H^\*(k[X],M)$, under the
$p$\nobreakdash-axis.

Convergence then implies that $H^0(A,M)\cong M^x$, that we have an
exact sequence
  \[\xymatrix@C-6pt{
  0 \ar[r]
    & H^1(A,M) \ar[r]
    & H^1(k[X],M) \ar[r]^-e
    & H^0(A,M) \ar[r]^-{d_2^{0,1}}
    & H^2(A,M) \ar[r]
    & 0
  }\]
with $e$ the edge morphism described in~\pref{p:edge} and that the
differentials provide isomorphisms
  \(
  d_2^{p,1}:H^p(A,M)\xrightarrow{\cong} H^{p+2}(A,M)
  \)
for each $p\geq1$.

According to~\pref{thm:E:hoch}, $d_2^{p,1}$ is given by Yoneda
multiplication by a class
  \[
  \zeta \in \Ext^2_{A^e}(\Tor_0^{k[X]}(A,A),\Tor_1^{k[X]}(A,A)) = \HH^2(A)
  \]
and, in fact, looking back at the proof of~\pref{thm:sOp=Ext}, we see
that $\zeta=d_2^{0,1}(1_A)$ with $1_A\in\HH^0(A)=A$.

On the other hand, it is easy to follow the recipe given
in~\pref{p:edge} in order to see that the edge morphism $e$ is induced
on $H^1(k[X],M)=M^x$ by the map $\hat e:M\to M_x$ given by
  \[
  \hat e(m) =
  \sum_{i=0}^N\sum_{\substack{s,t\geq0\\s+t+1=i}}\alpha_ix^smx^t.
  \]

\paragraph\label{p:pre:cup} In particular, if $M=A$ is the regular
$A$\nobreakdash-bimodule, of course we have that $A_x=A^x=A$, $e:A\to
A$ is just multiplication by $f'$, and we see that $\HH^0(A)\cong A$,
$\HH^1(A)\cong A_{f'}\denotes\ker e$, $\HH^2(A)\cong A/(f')$ and
that multiplication by $\zeta=d_2^{0,1}(1)\in\HH^2(A)$ gives an
epimorphism $\HH^0(A)\to\HH^2(A)$ and, for each $p\geq1$, an
isomorphism $\HH^p(A)\to\HH^{p+2}(A)$. 

\paragraph Under our isomorphisms, $\zeta$ corresponds to the class of
$1$ in $\HH^2(A)\cong A/(f')$. If we put $d=\gcd(f,f')$ and let $q$ be
such that $f=qd$, there is an isomorphism $A_{f'}\cong(q)/(f)$. In
particular $\HH^1(A)$ is a cyclic $A$\nobreakdash-module generated by
an element $\tau$, corresponding to the class of $q$ under this
isomorphism, with annihilator $(d)\lhd A$. 

\paragraph It is quite clear that the action of $\HH^0(A)=A$ on
$\HH^\*(A)$ corresponds, under our isomorphisms, to the obvious
structure of $A$\nobreakdash-modules on $A_{f'}$ and $A/(f')$. Since
we understand multiplication by $\zeta$, to describe the
multiplicative structure on $\HH^\*(A)$, we need only concentrate on
computing $\tau^2$.

\paragraph\label{p:obs:not2} Assume for a moment that $2$ is
invertible in $k$. The graded commutativity of $\HH^\*(A)$ and the
fact that $\abs\tau=1$ immediately imply, then, that $\tau^2=0$ and we
see that we have in this situation an isomorphism $\HH^\*(A)\cong
k[x,\tau,\zeta]/(f,\tau d,\zeta f')$ of graded commutative algebras.

\paragraph In the general case, we have the following theorem:

\begin{Theorem*}\label{thm:monogenic} Let $f=\sum_{i=0}^N\alpha_i
X^i\in k[X]$ be a monic polynomial of degree $N$ and consider the
$k$\nobreakdash-algebra $A=k[X]/(f)$. Let $d=\gcd(f,f')$, let $q\in
k[X]$ be such that $f=qd$ and put
  \[
  u= q^2 \sum_{i=0}^N \alpha_i\,\frac{i(i-1)}2\, X^{i-2}.
  \]
Then there is an isomorphism of graded commutative algebras
  \[
  \HH^\*(A)\cong k[x,\tau,\zeta]/(f(x),\tau d(x),\zeta f'(x),\tau^2-u(x)\zeta),
  \]
where the generators in the right hand side have degrees $\abs x=0$,
$\abs\tau=1$ and $\abs\zeta=2$.
\end{Theorem*}

\begin{proof}
At this point, we need only show that $\tau^2=u(x)\zeta$. We will
resort to a direct computation: we have not been able to find a more
conceptual argument in the spirit of those used above to handle this.

There is a commutative diagram of $A$\nobreakdash-bimodule morphisms
  \begin{equation}\label{eq:comp}
  \xymatrix@C-3pt{
  A\ot A \ar[r]^-{d_3} \ar[d]^-{s_3}
    & A\ot A \ar[r]^-{d_2} \ar@<0.5ex>[d]^-{s_2}
    & A\ot A \ar[r]^-{d_1} \ar@<0.5ex>[d]^-{s_1}
    & A\ot A \ar@{->>}[r]^-m \ar@<0.5ex>[d]^-{s_0}
    & A \ar@{=}[d]
    \\
  A\ot A^{\ot 3}\ot A \ar[r]^-{d'_3} 
    & A\ot A^{\ot 2}\ot A \ar[r]^-{d'_2} \ar@<0.5ex>[u]^-{r_2}
    & A\ot A\ot A \ar[r]^-{d'_1} \ar@<0.5ex>[u]^-{r_1}
    & A\ot A \ar@{->>}[r]^-m \ar@<0.5ex>[u]^-{r_0}
    & A
  }\end{equation}
where $m:A\ot A\to A$ is the multiplication map, the maps $d'_p$ are
the Hochschild boundary maps, $s_0=r_0=\id$ and
  \begin{gather*}
  \begin{aligned}
  &d_p(1\ot 1) = 1x\ot 1 - 1\ot x1 
        &&\qquad\text{if $p$ is odd,}\\
  &d_p(1\ot 1) = \sum_{i=0}^N\sum_{\substack{s,t\geq0\\s+t+1=i}}
                 \alpha_i\; 1x^s\ot x^t1  
       &&\qquad\text{if $p$ is even,}
  \end{aligned}
  \displaybreak[0] \\
  s_1(1\ot 1) = 1\ot x\ot 1,\\
  s_2(1\ot 1) = \sum_{i=0}^N\sum_{\substack{s,t\geq0\\s+t+1=i}}
        \alpha_i\; 1\ot x^s\ot x\ot x^t, \\
  s_3(1\ot 1) = \sum_{i=0}^N\sum_{\substack{s,t\geq0\\s+t+1=i}}
        \alpha_i\; 1\ot x\ot x^s\ot x\ot x^t,
        \displaybreak[0] \\
  r_1 (1\ot x^i\ot 1) = \sum_{\substack{s,t\geq0\\s+t+1=i}}
        x^s\ot x^t, \\
  r_2 (1\ot x^i\ot x^j\ot 1) = - 1\ot\quot_f(X^i+X^j).
  \end{gather*}
Here we are using functions $\quot_f,\,\rem_f:k[X]\to k[X]$ defined so
that, for all $h\in k[X]$, $h=\quot_f(h)f+\rem_f(h)$ and either
$\rem_f(h)=0$ or $\deg\rem_f(h)<\deg f$.

The rows in~\eqref{eq:comp} are exact and, in fact, they are the
beginnings of two projective resolutions of $A$ as an
$A$\nobreakdash-bimodule: the lower row comes from the usual
Hochschild resolution $A^{\ot(\*+2)}$ of $A$ and the upper row comes
from the well-known $2$\nobreakdash-periodic resolution $P_\*$ of $A$
constructed in~\cite{BACH}. The vertical maps are the first components
of a comparison of resolutions. The complete picture can be found
in~\cite{Holm}, but we will not make use of it.

It is clear that $\tau\in\HH^1(A)$ is the class of the unique
derivation $t:A\to A$ such that $t(x)=q$. This implies that $\tau$ can
be seen as the class of the $1$\nobreakdash-cocycle $t\in\hom(A,A)$ of
the complex $\hom(A^{\ot\*},A)\cong\hom_{A^e}(A^{\ot(\*+2)},A)$
constructed by applying $\hom_{A^e}(\place,A)$ to the Hochschild
resolution of $A$.

Now, when one sees $\HH^\*(A)$ as the cohomology of this complex,
products can be computed using the cup product~$\smile$ introduced
in~\cite{Gerstenhaber}. This means that $\tau^2$ is the class of
$t\smile t:A\ot A\to A$ in $\HH^2(A)$, where
  \[
  (t\smile t)(a\ot b)=t(a)t(b) = a'b'q^2,
  \]
for all $a,b\in A$; here derivatives are taken on arbitrary
representatives for elements of $A$ in $k[X]$.

Pulling back the $2$\nobreakdash-cocycle $t\smile t$ from
$\hom_{A^e}(A^{\ot(\*+2)},A)$ to the top row in~\eqref{eq:comp} along
the given comparison morphisms immediately shows that $\tau^2$ is
represented by the unique $2$\nobreakdash-cocycle of the complex
$\hom_{A^e}(P_\*,A)$ which maps $1\ot1\in P_2=A\ot A$ to
  \[
  \sum_{i=0}^N\sum_{\substack{s,t\geq0\\u+v+1=i}}
        \alpha_i\;t(x^u)t(x)x^v
  =
  \sum_{i=0}^N\sum_{\substack{s,t\geq0\\u+v+1=i}}
        \alpha_i\,u\, x^{i-2}q^2
  =
  \sum_{i=0}^N
        \alpha_i\,\frac{i(i-1)}2\, x^{i-2}q^2.
  \]
Since the class $\zeta$ is represented in the complex
$\hom_{A^e}(P_\*,A)$ by the $2$\nobreakdash-cocycle $z:A\ot A\to A$
such that $z(1\ot 1)=1$, we see that $t\smile t=u(x)z$ and, then,
that $\tau^2=u(x)\zeta$, as stated in~\pref{thm:monogenic}. This
completes the proof of that theorem.
\end{proof}

\paragraph If $2$ is invertible in $k$, then the polynomial $u$
appearing in~\pref{thm:monogenic} is simply $\frac12q^2f''$ and it is
an easy exercise to show that this is zero in $A/(f')$. This
corresponds, of course, to the observation made in~\pref{p:obs:not2}.

\paragraph Before passing on to other matters, we take the opportunity
of computing the rest of the ``cohomology structure'' of our algebra
$A$ in the sense used in~\cite{Gerstenhaber}. We will use the
conventions and notations of~\cite{GerstenhaberSchack1}*{Section~4},
which are, by now, rather standard; in particular, we use the
composition products $\circ$ and $\circ_i$ for cocycles on the
Hochschild resolution.

\begin{Theorem*}\label{thm:monogenic:lie}
In the situation of theorem~\pref{thm:monogenic}, let
  \[
  w = \sum_{i=0}^N\sum_{\substack{s,t\geq0\\s+t+1=i}}
         \alpha_i\,\quot_f((s+1)x^sq)x^t.
  \]
Then the Gerstenhaber bracket on~$\HH^\*(A)$ is
completely determined by the relations
  \begin{gather*}
  [\tau,x] = q, \\
  [\zeta,\tau] = w\zeta, \\
  [x,x] = [\tau,\tau] = [\tau,\zeta] = [x,\zeta] = 0.
  \end{gather*}
\end{Theorem*}
\unskip

\begin{proof}
The Gerstenhaber bracket is graded on $\HH^\*(A)[1]$, so that it is
clear that $[x,x]\in\HH^{-1}(A)$ and $[\tau,\tau]$ are zero. The other
four relations will be established by computation. We use the maps
$r_\*$ and $s_\*$ to go from cocycles on the top row
of~\eqref{eq:comp} to cocycles on the bottom row and back,
respectively, and identify $\hom_{A^e}(A^{\ot(\*+2)},A)$ with
$\hom(A^{\ot\*},A)$ as usual. Also, we will write 
  \[
  \int \phi(i,s,t) 
  \denotes 
  \sum_{i=0}^N\sum_{\substack{s,t\geq0\\s+t+1=i}}
  \phi(i,s,t)
  \]
for functions $\phi$ defined on non-negative integers.

We have that
  \begin{align*}
  s_0^*\big([r_0^*(x),r_1^*(\tau)]\big)(1\ot 1)
      &= [r_0^*(x),r_1^*(\tau)](1) \\
      &= \big(r_1^*(\tau)\circ(r_0^*(x)\big)(1) \\
      &= r_1^*(\tau)\big(r_0^*(x)(1)\big) \\
      &= r_1^*(\tau)(x) \\
      &= q,
  \end{align*}
so $[x,\tau]=q$. Similarly, we see that
  \begin{align*}
  s_1^*\big([r_0^*(x),r_2^*(\zeta)]\big)(1\ot 1)
      &= [r_0^*(x),r_2^*(\zeta)](x) \\
      &= \big(r_2^*(\zeta)\circ r_0^*(x)\big)(x) \\
      &= r_2^*(\zeta)\big(r_0^*(x),x)-r_0^*(\zeta)(x,r_0^*(x)) \\
      &= 0,
  \end{align*}
because $r_2^*(\zeta):A^{\ot2}\to A$ is symmetric. This tells us
that $[x,\zeta]=0$. 

Using symmetry again, we compute
  \begin{align*}
  s_3^*&\big([r_2^*(\zeta),r_2^*(\zeta)]\big)(1\ot 1)
      = \int
          \alpha_i\,[r_2^*(\zeta),r_2^*(\zeta)](x\ot x^s\ot x)\,x^t 
          \displaybreak[0] \\
    & = 2 \int
          \alpha_i\,(r_2^*(\zeta)\circ r_2^*(\zeta))(x\ot x^s\ot x)\,x^t 
          \displaybreak[0] \\
    & = 2 \int
            \alpha_i\,\Big(
                r_2^*(\zeta)\big(r_2^*(\zeta)(x\ot x^s)\ot x\big) 
              \mathop- r_2^*(\zeta)\big(x\ot r_2^*(\zeta)(x\ot x^s)\big)
                    \Big)\,x^t
          \\
    & = 0,
  \end{align*}
so we see that $[\zeta,\zeta]=0$.

Finally, 
  \begin{align*}
  s_2^*\big([r_1^*(\tau),r_2^*(\zeta)]\big)(1\ot 1)
     &= \int
         \alpha_i\,[r_1^*(\tau),r_2^*(\zeta)](x^s\ot x) x^t 
       \\
   & = \int
         \alpha_i\,\left(\quot_f((s+1)x^sq)-\quot_f(x^{s+1})'q\right)x^t
\intertext{and, since $\quot_f(x^{s+1})'=0$ if $0\leq s<N$, this is}
   & = \int
         \alpha_i\,\quot_f((s+1)x^sq)x^t.
  \end{align*}
We have thus verified all the relations claimed in the theorem.
\end{proof}

\section{Variations}
\label{sect:variations}

The computation done in section~\pref{sect:monogenic} was successful
because of the many favorable traits of the situation under
consideration. It turns out, though, that a similar line of reasoning
can be applied in various other less favorable contexts in order to
obtain useful information on cohomology. We collect here a few
examples.

\begin{Proposition}\label{p:flat:ideal}
Let $A$ be an algebra and $I\lhd A$ an ideal which is flat as an
$A$\nobreakdash-module on the left or on the right and put $B=A/I$.
Then $H^0(B,\place)\cong H^0(A,\place)$ and there is a long exact
sequence
  \begin{gather*}
  \xymatrix@C-12pt{
    0  \ar[r]
      & H^1(B,M) \ar[r]
      & H^1(A,M) \ar[r]
      & \Ext_{B^e}^0(I/I^2,M) \ar[r]^-{\smile\zeta}
      & H^2(B,M) \ar[r]
      & \cdots
  } \\
  \mskip-20mu
  \xymatrix@C-12pt{
    \cdots  \ar[r]
      & H^p(B,M) \ar[r]
      & H^p(A,M) \ar[r]
      & \Ext_{B^e}^{p-1}(I/I^2,M) \ar[r]^-{\smile\zeta}
      & H^{p+1}(B,M) \ar[r]
      & \cdots
  }
  \end{gather*}
functorial on $B$\nobreakdash-bimodules $M$. Here $\zeta\in
H^2(B,I/I^2)$ is the class of the singular extension 
  \[\xymatrix{
  0 \ar[r]
    & I/I^2 \ar[r]
    & A/I^2 \ar[r]
    & B \ar[r]
    & 0
  }\]
In particular, restriction of scalars induces functorial isomorphisms
  \[
  H^p(B,\place) \to H^p(A,\place)
  \]
on $B$\nobreakdash-bimodules for $p>\pdim_{A^e}I/I^2+1$.
\end{Proposition}

\begin{proof}
Looking at the long exact sequence for $\Tor_\*^A(B,\place)$
corresponding to 
  \[\xymatrix{
  0\ar[r] &  I\ar[r] &  A\ar[r] &  B\ar[r] &  0
  }
  \]
we see that $\Tor_p^A(B,B)\cong\Tor_{p-1}^A(B,I)=0$ for $p>1$. Using
this together with the convergence of the spectral sequence
in~\pref{thm:E:hoch} we see that the long exact sequence in the
statement exists. All other claims follow at once.
\end{proof}

\paragraph\label{p:normal} Let $A$ be an algebra, let $x\in A$ be
\definiendum{normal} (so that $Ax=xA$) and assume moreover that $x$ is
not a divisor of zero. Normality implies that the left ideal ${I=Ax}$ is
actually a bilateral ideal and we can consider the quotient algebra
${B=A/I}$. The map $r_x:A\to I$ given by right multiplication by~$x$
is an isomorphism of left $A$-modules\dash---in particular, $I$~is
flat and we can apply~\pref{p:flat:ideal} to this situation. 

One can see at once that the hypothesis on~$x$ implies that there
exists a unique automorphism $\alpha\in\Aut_\Alg(A)$ such that
$ax=x\alpha(a)$ for all $a\in A$; notice that $\alpha=\id_A$ iff
$x$~is central. Moreover, $\alpha(I)=I$ so $\alpha$~induces an
automorphism of~$B$, which we denote~$\alpha$ as well. If~$M$ is a
$B$-bimodule, we write $M_\alpha$ the $B$-bimodule which coincides
with $M$ as a left $B$-module and whose right action is that of~$M$
`twisted' by~$\alpha$, so that
  \[
  m\cdot b=m\alpha(b), \qquad\forall m\in M_\alpha,\,b\in B.
  \]

Clearly, $r_x(I)=I^2$ and in fact $r_x$~induces an isomorphism of
$B$-bimodules $I/I^2\cong B_\alpha$. Now, it is easy to see that for
each $B$-bimodule~$M$, there is a natural isomorphism
$\Ext_{B^e}^\*(B_\alpha,M) \cong\Ext_{B^e}^\*(B,M_{\alpha^{-1}})
=H^\*(B,M_{\alpha^{-1}})$; indeed, this follows from the fact that the
functor $(\place)\ot_B B_\alpha:\bMod BB\to\bMod BB$ is an equivalence
which maps~$B$ to~$B_\alpha$ and $M_{\alpha^{-1}}$ to~$M$. Using this,
the long exact sequence of~\pref{p:flat:ideal} becomes
  \begin{gather}\label{eq:balpha:lec}
  \hskip-3pt 
  \xymatrix@C-12pt{
    0  \ar[r]
      & H^1(B,M) \ar[r]
      & H^1(A,M) \ar[r]
      & H^0(B,M_{\alpha^{-1}}) \ar[r]^-{\smile\zeta}
      & H^2(B,M) \ar[r]
      & \cdots
  } \\
  \hskip-3pt 
  \xymatrix@C-12pt{
    \cdots  \ar[r]
      & H^p(B,M) \ar[r]
      & H^p(A,M) \ar[r]
      & H^{p-1}(B,M_{\alpha^{-1}}) \ar[r]^-{\smile\zeta}
      & H^{p+1}(B,M) \ar[r]
      & \cdots
  } \notag
  \end{gather}
The following proposition records an interesting special case which
occurs when $\pdim_{A^e}A\leq 1$\dash---for example, when $A$ is finite
dimensional and hereditary.

\begin{Proposition}\label{p:balpha}
Let $A$ be an algebra such that $\pdim_{A^e}A\leq1$. Let $x\in A$ be a
normal non-zero divisor and let $\alpha\in\Aut_\Alg(A)$ be such that
$ax=x\alpha(a)$ for all $a\in A$. Put $I=(x)$ and $B=A/I$. Then for
each $B$-bimodule~$M$ and each $p\geq0$ there is an exact sequence
  \begin{multline*}
  \xymatrix@C-5pt{
    0  \ar[r]
      & H^{2p+1}(B,M) \ar[r]
      & H^1(A,M_{\alpha^{-p}}) \ar[r]^-e
      & {}
  } \\
  \xymatrix@C-5pt{
    {} \ar[r]
      & H^0(A,M_{\alpha^{-p-1}}) \ar[r]
      & H^{2p+2}(B,M) \ar[r]
      & 0
  } \qedhere
  \end{multline*}
\end{Proposition}

\noindent We thus see that under the stated conditions, computation of
cohomology is reduced to the consideration of what happens in low
degrees.

\begin{proof}
As $\pdim_{A^e}A\leq1$, the long exact sequence~\eqref{eq:balpha:lec}
collapses into an exact sequence
  \[\xymatrix@C-5pt{
    0  \ar[r]
      & H^1(B,M) \ar[r]
      & H^1(A,M) \ar[r]
      & H^0(B,M_{\alpha^{-1}}) \ar[r]
      & H^2(B,M) \ar[r]
      & 0
  }\]
and natural isomorphisms $H^p(B,M_{\alpha^{-1}})\cong H^{p+2}(B,M)$
for all $p\geq1$. Iterating these isomorphisms, using the exact
sequence and taking into account the the isomorphism
$H^0(A,\place)\cong H^0(B,\place)$,  we obtain the sequences referred
to in the proposition.
\end{proof}

\paragraph The proof of the proposition also shows the following:

\begin{Corollary*}
With the notations of~\pref{p:balpha}, assume that $\alpha^n=\id_A$.
Then there is a class $\xi\in\HH^{2n}(B)$ such that
$\smile\xi:H^{p}(B,M)\to H^{p+2n}(B,M)$ is an isomorphism for all
$p\geq1$ and all $B$-bimodules~$M$. 
In particular, if $\dim_kB<\infty$, then $\HH^\*(B)$ is a
finitely generated algebra and $H^\*(B,M)$ is a finitely generated
$\HH^*(B)$-module for all finite dimensional $B$-bimodules~$M$.
\end{Corollary*}

\begin{proof}
The class $\zeta\in H^2(B,I/I^2)$ of the extension
  \[\xymatrix{
  0 \ar[r]
    & I/I^2 \ar[r]
    & A/I^2 \ar[r]
    & B \ar[r]
    & 0
  }\]
is such that $\smile\zeta: H^p(B,M)\to H^{p+2}(B,M_\alpha)$ is an
isomorphism for all $B$-bimodules~$M$ and all $p\geq1$. Now the $n$-th
Yoneda power of~$\zeta$ can be seen as a class $\xi=\zeta^n\in
H^2(B,(I/I^2)^{\ot_Bn})=H^2(B,B_{\alpha^n})=\HH^{2n}(B)$ and, of
course, the map $\smile\xi:H^p(B,M)\to H^{p+2n}(B,M)$ is an
isomorphism for all $p\geq0$. This proves the first claim in the
statement.

Assume now that $\dim_kB<\infty$ and let $M$ be a finite dimensional
$B$-bimodule. What we have so far implies that $H^\*(B,M)$ is
generated as a $\HH^\*(B)$\nbd-module by $\bigoplus_{p=0}^{2n-1}H^p(B,M)$,
which is a finite dimensional vector space. It follows at once that
$H^\*(B,M)$~is a finitely generated module. Similarly, $\HH^\*(B)$~is
generated as an algebra by~$\bigoplus_{p=0}^{2n-1}\HH^p(B)$ together
with~$\xi$, so it is itself finitely generated.
\end{proof}

\paragraph In some cases, we can obtain more precise information about
the multiplicative structure on Hochschild cohomology in the situation
of~\pref{p:normal}. We consider here only a very simple instance.

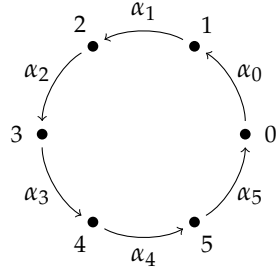
\begin{figure}[t]
\begin{tikzpicture}[scale=1.35]
  \def\last{-60}
  \foreach \i / \j / \b in {0/5/0, 1/0/60, 2/1/120, 3/2/180, 4/3/240, 5/4/300}
    {
      \fill (\b:1) circle (1.5pt) ;
      \draw (\b:1.25) node {\small$\i$} ;
      \draw[->, shorten <=5pt, shorten >=5pt] 
        (\last:1) arc (\last:\b:1) ;
      \path (\last:1.4) -- node {$\alpha_\j$} (\b:1.4) ;
      \global\let\last\b
    }
\end{tikzpicture}
\caption{The $6$-crown.}
\label{fig:taft}
\end{figure}

Let $n\in\NN$ and let $Q$ be the quiver with vertex set $Q_0=\ZZ_n$
and arrow set $Q_1=\{\alpha_i\}_{i\in\ZZ_n}$ such that $\alpha_i$ starts
at the vertex~$i$ and ends at the vertex~$i+1$. This is sometimes
called the \definiendum{$n$-crown quiver}. See figure~\ref{fig:taft}
for a drawing of~$Q$ when~$n=6$.

Let $kQ$ be the path algebra on $Q$ and put $t=\sum_{\alpha\in
Q_1}\alpha\in kQ$. If $f\in k[X]$ is a non-zero polynomial, one easily
sees that $u=f(t)\in kQ$ is a regular normal element and if we put
$I=(u)$ and $B=kQ/I$, the observations in~\pref{p:normal} provide a
long exact sequence relating $\HH^\*(B)$ and $H^\*(kQ,B)$. 

We take $f=X^l$ with $l=nm-1$ for some $m\geq1$, so that $B$ is one of
the symmetric truncated cycle algebras considered in~\cite{BLM} or
the special Brauer tree algebras studied in~\cite{Holm:Brauer}.
Additionally, we assume $n\geq2$: the case in which $n=1$ was treated
in the previous section.

The choice of~$l$ implies that $\alpha=\id_A$. Using this and that
$\pdim_{(kQ)^e}kQ=1$, we see that the long exact
sequence~\eqref{eq:balpha:lec} degenerates, as in the proof
of~\pref{p:balpha}, into, on one hand, an exact sequence
  \begin{equation}\label{eq:cycle:suc}
  \xymatrix@C-5pt{
  0 \ar[r]
    & \HH^1(B) \ar[r]
    & H^1(kQ,B) \ar[r]^-e
    & \HH^0(B) \ar[r]^-{\smile\zeta}
    & \HH^2(B) \ar[r]
    & 0
  }\end{equation}
with $e$ the edge morphism from the limit to the `base' described
in~\pref{p:edge} and $\zeta\in\HH^2(B)$ the class of the extension
of algebras
  \begin{equation}\label{eq:cycle:ext}
  \xymatrix@C-5pt{
  0 \ar[r]
    & B\cong I/I^2 \ar[r]
    & kQ/I^2 \ar[r]
    & B \ar[r]
    & 0
  }\end{equation}
and, on the other hand, isomorphisms $\HH^p(B)\to\HH^{p+2}(B)$ for
$p\geq1$, all given by cup-product with~$\zeta$. 

\begin{Proposition*}
There is an isomorphism of algebras $\HH^\ev(B)\cong\Z(B)[\zeta]$,
there is an isomorphism of $\HH^\ev(B)$-modules $\HH^\odd(B)\cong
H^1(kQ,B)[\zeta]$ and the product of two elements of $\HH^\odd(B)$ is
zero. This information, together with the usual $\Z(B)$-bimodule
structure on $H^1(kQ,B)$, completely determines the cohomology algebra
  \[
  \HH^\*(B)\cong\Z(B)[\zeta]\oplus H^1(kQ,B)[\zeta].
  \]
\par\removelastskip
\end{Proposition*}

\begin{proof}
Let $E\subset kQ$ be the subalgebra generated by the vertices. One
easily checks that, because $n\geq2$, every $E$-linear derivation
$kQ\to B$ into the $kQ$-bimodule~$B$ induces an $E$-linear derivation
$B\to B$ simply by composing with the projection~$\phi$. This means
that the map $\HH^1(B)\to H^1(kQ,B)$ in~\eqref{eq:cycle:suc} is
surjective and, of course, it follows that $e=0$ and that
$\smile\zeta:\HH^0(B)\to\HH^2(B)$ is an isomorphism. We can conclude
that $\HH^\ev(B)\cong\Z(B)[\zeta]$. We already know how $\zeta$ acts
on $\HH^\odd(B)$ and the center $\Z(B)=\HH^0(B)$ acts on it
canonically, so we need only describe the restriction of the
cup-product to $\HH^1(B)\ot\HH^1(B)$ in order to complete the proof.

Let us write $\gamma_i^j$ the only path in $Q$ of length $j$ which
starts in the vertex $i\in\ZZ_n$. Then
$\{\gamma_i^j:i\in\ZZ_n,\;j\geq0\}$ is a basis of $kQ$ and
$\{\gamma_i^j:i\in\ZZ_n,\;0\leq j<l\}$ is a basis of $B$. Let $kQ_1$
and $R$ be the sub-$E^e$-modules of $kQ$ spanned by $Q_1$ and
$\{\gamma_i^l:i\in\ZZ_n\}$, respectively. We consider the following
commutative diagram
  \[\xymatrix{
  B\ot_ER\ot_EB \ar[r]^-{d_2} \ar[d]^-{s_2}
    & B\ot_E kQ_1\ot_EB \ar[r]^-{d_1} \ar@<0.5ex>[d]^-{s_1}
    & B\ot_E B \ar@{->>}[r]^-\mu \ar@{=}[d]
    & B \ar@{=}[d]
    \\
  B^{\ot_E4} \ar[r]
    & B^{\ot_E3} \ar[r] \ar@<0.5ex>[u]^-{r_1}
    & B\ot_E B \ar@{->>}[r]^-\mu 
    & B 
  }\]
in which the bottom row is the standard Hochschild resolution of $B$
taken over~$E$, the map~$\mu$ is the multiplication, $s_1$ is the
inclusion and 
  \begin{gather*}
  d_1(1\ot\alpha\ot 1) = \alpha\ot 1-1\ot\alpha,
    \\
  d_2(1\ot\gamma_i^l\ot 1) 
    = \sum_{s=0}^{l-1} 
        \gamma_{i+s+1}^{l-s-1}\ot\alpha_{i+s}\ot\gamma_i^s,
    \\
  r_1(1\ot\gamma_i^j\ot1) 
    = \sum_{s+t+1=j} \gamma_{i+t+1}^{s}\ot\alpha_{i+t}\ot\gamma_i^t,
    \\
  s_2(1\ot\gamma_i^l\ot 1)
    = \sum_{s+t+1=l}
        1\ot\gamma_{i+t+1}^{s}\ot\alpha_{i+t}\ot\gamma_i^t,
  \end{gather*}
The top row is the beginning of the Bardzell resolution for $B$,
\cf~\cite{Bardzell}.

Let us take now classes $\phi$, $\psi\in\HH^1(B)$ and $1$-cocycles
$f$, $g:B\ot_EkQ_1\ot_EB\to B$ defined on the Bardzell complex such
that $f$ and $g$ represent $\phi$ and $\psi$, respectively. We
identify $f$ and $g$ canonically to elements of $\hom_{E^e}(Q_1,B)$,
which can itself be seen as the vector space of $E$-linear derivations
$kQ\to B$. Then $s_2^*(r_1^*(f)\smile r_1^*(g))$ is a representative
for $\phi\smile\psi$. Since
  \begin{multline*}
  s_2^*\big(r_1^*(f)\smile r_1^*(g)\big)(1\ot\gamma_i^l\ot 1)
    \\
    = \sum_{s+1+t+1+v=l}
        \gamma_{i+1+t+1+v}^s
        f(\alpha_{i+t+1+v})
        \gamma_{i+1+v}^t
        g(\alpha_{i+v})
        \gamma_i^v
  \end{multline*}
and every element of $\hom_{E^e}(Q_1,B)$ maps arrows to linear
combinations of paths of length at least one (because $n\geq2$) we see
that $s_2^*\big(r_1^*(f)\smile r_1^*(g)\big)=0$, so in particular,
$\phi\smile\psi=0$ in~$\HH^2(B)$. This implies that the cup product
vanishes on~$\HH^\odd(B)$.
\end{proof}


\section{The cohomology of nice quotients}
\label{sect:nice}

\paragraph It is well know that a morphism of algebras $\phi:A\to B$
is an epimorphism in the category $\Alg$ of algebras iff extension of
scalars along $\phi$ is a full and faithful functor $\phi^*:\lMod
B\to\lMod A$; in particular, a surjection is an epimorphism.

Following~\cite{GeigleLenzig}, one says that $\phi$ is a
\definiendum{homological epimorphism} if the corresponding functor on
bounded derived categories $D^b(\phi^*):D^b(\lMod B)\to D^b(\lMod A)$
is full and faithful. An easy induction on the length of complexes
shows that this is equivalent to the condition that extension of
scalars induce an isomorphism
$\Ext_B^\*(\place,\place)\to\Ext_A^\*(\place,\place)$ of bifunctors of
$B$\nobreakdash-modules.

\paragraph The kernel of a homological epimorphism is called
a~\definiendum{homological ideal}. 

Homological ideals were first considered by Maurice Auslander~\etal{}
in~\cite{APT} under the name of~\definiendum{strong idempotent
ideals}: indeed, it turns out that ideals of this form are idempotent,
\cf~\pref{p:hepi-is-idemp}.

\paragraph The characterization of homological epimorphisms given in
the following proposition comes from~\cite{GeigleLenzig}
and~\cite{delaPenhaXi}*{Proposition~2.3}, though the fact that
$\text{\ref{hepi:4}}\Rightarrow\text{\ref{hepi:1}}$ was not noted
there. The final statement about Hochschild cohomology
strengthens~\cite{delaPenhaXi}*{Proposition~3.1}; in particular,
remark that we are not restricting ourselves to finite dimensional
algebras.

\begin{Proposition*}\label{prop:hepi}
Let $\phi:A\to B$ be an morphism of algebras such that $B\ot_AB\cong
B$ as $B$\nobreakdash-bimodules. The following statements are equivalent:
\begin{enumerate}[label=(\alph*),ref=\emph{(\alph*)}]
\item\label{hepi:1} $\phi:A\to B$ is a homological epimorphism;
\item\label{hepi:2} $\Tor_+^A(B,M)=0$ for all $M\in\lMod B$;
\item\label{hepi:3} $\Tor_+^A(B,B)=0$;
\item\label{hepi:4} $\phi^e:A^e\to B^e$ is a homological epimorphism.
\end{enumerate}
When they are satisfied, there is a natural isomorphism
  \begin{equation}\label{prop:hepi:iso}
  H^\*(B,\place)\to H^\*(A,\place)
  \end{equation}
of functors of $B$\nobreakdash-bimodules. Furthermore, the
specialization of~\eqref{prop:hepi:iso} to the
$B$\nobreakdash-bimodule~$B$ is an isomorphism of algebras
$\HH^\*(B)\to H^\*(A,B)$.
\end{Proposition*}

\begin{proof}
\def\imp#1=>#2{$\text{\ref{hepi:#1}}\Rightarrow\text{\ref{hepi:#2}}$}
We prove \imp1=>2. Fix $M\in\lMod B$ and consider the functorial
spectral sequence with
$E_2^{p,q}(\place)\cong\Ext_B^p(\Tor_q^A(B,M),\place)\Rightarrow
\Ext_A^\*(M,\place)$ constructed in~\pref{prop:change}. The edge
morphisms from the fiber to the limit in this sequence are then maps
$e:E_2^{\*,0}(\place)\cong\Ext_B^\*(M,\place)\to\Ext_A^\*(M,\place)$
which we are assuming to be isomorphisms.

Assume that $r\geq0$ and that we know that $\Tor_q^A(B,M)=0$ if
$0<q<r$. Then $E_{r+1}^{p,q}=0$ if $0<q<r$,
$E_{r+1}^{p,q}=E_2^{p,q}$ if $(p,q)\in\{(r,0),(0,r),(r+1,0)\}$, and 
convergence implies that we have an exact sequence
 \[\xymatrix@C-9pt{
 E_{r+1}^{r,0}(\place) \ar[r]^-e
   & \Ext_A^r(M,\place) \ar[r]
   & E_{r+1}^{0,r}(\place) \ar[r]^-{d_{r+1}^{0,r}}
   & E_{r+1}^{r+1,0}(\place) \ar[r]^-e
   & \Ext_A^{r+1}(M,\place)
 }\]
This tells us that $E_2^{0,r}(\place)=\hom_B(\Tor_r^A(B,M),\place)$
vanishes identically on $\lMod B$ and of course allows us to conclude
that $\Tor_r^A(B,M)=0$. 
This argument clearly gives us~\ref{hepi:2} by induction.

The implication \imp2=>3 is immediate. To show that \imp3=>1, we
remark that~\ref{hepi:3} implies that the spectral sequence
constructed in~\pref{thm:E:hoch} degenerates, so that the
edge morphisms $\Ext_{B^e}^\*(B,\place)\to\Ext_{A^e}^\*(A,\place)$ are
isomorphisms on $\bMod BB$; \ref{hepi:1}~now follows from the easily
established fact that for an algebra $\Lambda$ there is a natural
isomorphism 
  \begin{equation}\label{eq:ext:as:hoch}
  \Ext_\Lambda^\*(\place,\place)\cong H^\*(\Lambda,\hom(\place,\place))
  \end{equation}
of bifunctors defined on $\lMod\Lambda$, which is compatible with
extension of scalars.

Recall that there is a map
$\top:\Tor^A_\*(B,B)\ot\Tor^{A^{\op}}_\*(B,B)\to\Tor^{A^e}_\*(B^e,B^e)$
which is an isomorphism,
\cf~\cite{CartanEilenberg}*{Theorem~\textsc{XI}.3.1}. From this and the
obvious existence of an isomorphism $\Tor_\*^{A^\op}(B,B)\cong\Tor_\*^A(B,B)$ we see
that~\ref{hepi:3} implies that $\phi^e:A^e\to B^e$ itself
satisfies~\ref{hepi:3}. Using what we have already proved, we
conclude that~\imp3=>4.

To show~\imp4=>1, assume $\phi^e$ is a homological epimorphism.
Applying the implication~\imp1=>3 to $\phi^e$, we get that
  \[
  \big(\Tor^A_\*(B,B)\ot\Tor^{A^{\op}}_\*(B,B)\big)_+
    \cong \Tor^{A^e}_+(B^e,B^e) 
    = 0,
  \]
so $\Tor^A_+(B,B)=0$. Now we can use the fact that~\imp2=>1 to
conclude that $\phi$ is a homological epimorphism.

Finally, the statement about Hochschild cohomology follows from the
fact that~\ref{hepi:3} implies that the spectral sequence
in~\pref{thm:E:hoch} degenerates and the statement about cup products
follows from~\pref{p:edge:cup}.
\end{proof}

\paragraph The following is an easy corollary of~\pref{prop:hepi}:

\begin{Corollary*}
Let $\phi:A\to B$ be a homological epimorphism. If $I=\ker\phi$, then
$\Ext_{A^e}^\*(I,\place)$ vanishes identically on $\bMod BB$. 
\end{Corollary*}

\begin{proof}
Let $M\in\bMod BB$ and
consider the long exact sequence of the functor
$\Ext_{A^e}^\*(\place,M)$ corresponding to $0\to I\to A\to B\to 0$.
The composition of the induced map
  \[
  \Ext_{A^e}^\*(B,M)\to\Ext_{A^e}^\*(A,M)
  \]
with the isomorphism $\Ext_{B^e}^\*(B,M)\to\Ext_{A^e}^\*(B,M)$
coincides with the isomorphism~\eqref{prop:hepi:iso}, so we see that
$\Ext_{A^e}^\*(I,M)=0$.
\end{proof}

\paragraph\label{p:pre-happel} Just as easily, we obtain the following
much more interesting result. This is what the first two parts
of~\cite{delaPenhaXi}*{Proposition~3.2} become when
taking~\pref{prop:hepi} into account. 

\begin{Corollary*}
Let $\phi:A\to B$ be a homological epimorphism and let $I=\ker\phi$.
There is a long exact sequence
  \begin{multline*}
  \xymatrix@C-10pt{
    0 \ar[r]
      & H^0(A,I) \ar[r]
      & \HH^0(A) \ar[r]
      & \HH^0(B) \ar[r]
      & \Ext_{A^e}^1(A,I) \ar[r]
      & \cdots
      } \\
  \xymatrix@C-10pt{
    \cdots \ar[r]
      & \Ext_{A^e}^p(A,I) \ar[r]
      & \HH^p(A) \ar[r]
      & \HH^p(B) \ar[r]
      & \Ext_{A^e}^{p+1}(A,I) \ar[r]
      & \cdots
      }
  \end{multline*}
The maps $\HH^p(A)\to\HH^p(B)$ appearing in this sequence 
can be collected into an algebra morphism $\HH^\*(A)\to\HH^\*(B)$.
\end{Corollary*}

\begin{proof}
This is just the long exact sequence for the
functor~$\hom_{A^e}(A,\place)$ corresponding to the short exact
sequence 
  \[\xymatrix{
  0\ar[r] &  I\ar[r] &  A\ar[r] &  B\ar[r] &  0
  }\]
up to the isomorphism $\HH^\*(B)\cong H^\*(A,B)$ coming
from~\eqref{prop:hepi:iso}.

The last statement is a direct consequence of the last statement
in~\pref{prop:hepi}.
\end{proof}

\paragraph\label{p:hepi-is-idemp} If $\phi:A\to B$ is an epimorphism
and $I=\ker\phi$, $\Tor_1^A(B,B)\cong I/I^2$, so~\pref{prop:hepi}
implies that a homological ideal is idempotent. We have the following
partial converse:

\begin{Proposition*}
Let $I\lhd A$ be flat as an $A$\nobreakdash-module on the left or on
the right and put $B=A/I$. If it is idempotent, then it is homological.
In general, $H^0(B,\place)\cong H^0(A,\place)$ and there is a long
exact sequence
  \begin{gather*}
  \xymatrix@C-12pt{
    0  \ar[r]
      & H^1(B,M) \ar[r]
      & H^1(A,M) \ar[r]
      & \Ext_{A^e}^0(I/I^2,M) \ar[r]^-{\smile\zeta}
      & H^2(B,M) \ar[r]
      & \cdots
  } \\
  \mskip-20mu
  \xymatrix@C-12pt{
    \cdots  \ar[r]
      & H^p(B,M) \ar[r]
      & H^p(A,M) \ar[r]
      & \Ext_{A^e}^{p-1}(I/I^2,M) \ar[r]^-{\smile\zeta}
      & H^{p+1}(B,M) \ar[r]
      & \cdots
  }
  \end{gather*}
functorial on $B$\nobreakdash-bimodules $M$. Here $\zeta\in
H^2(B,I/I^2)$ is the class of the singular extension 
  \[\xymatrix{
  0 \ar[r]
    & I/I^2 \ar[r]
    & A/I^2 \ar[r]
    & B \ar[r]
    & 0
  }\]
In particular, restriction of scalars induces functorial isomorphisms
  \[
  H^p(B,\place) \to H^p(A,\place)
  \]
on $B$\nobreakdash-bimodules for $p>\pdim_{A^e}I/I^2+1$.
\end{Proposition*}

\begin{proof}
Looking at the long exact sequence for $\Tor_\*^A(B,\place)$
corresponding to 
  \[\xymatrix{
  0\ar[r] &  I\ar[r] &  A\ar[r] &  B\ar[r] &  0
  }
  \]
we see that $\Tor_p^A(B,B)\cong\Tor_{p-1}^A(B,I)=0$ for $p>1$. Using
this together with the convergence of the spectral sequence
in~\pref{thm:E:hoch} we see that the long exact sequence in the
statement exists. All other claims follow at once.
\end{proof}

\paragraph We now describe a nice example where one can
see~\pref{p:hepi-is-idemp} in nature. Let $X$ be a finite poset and
let $Y\subset X$ be an order ideal, that is, a subset such that
  \[
  x\in X,\; y\in Y,\;x\leq y
  \implies
  x\in Y.
  \]
Let $A=kX$ and $kY$ be the incidence algebras of $X$ and
$Y$\dash---recall that $kX$, for example, can be seen as the quotient
of the path algebra on the quiver given by the Hasse diagram of $X$
divided by the ideal of all commutativity relations. For simplicity,
we identify $X$ with (the vertex set of) its Hasse diagram.

Let us put $e=\sum_{x\in X\setminus Y}x$. This is an idempotent in
$A$, so $I_Y=AeA$ is an idempotent ideal. Now, $I_Y$ is clearly
linearly spanned by all paths in $X$ which go through a vertex in
$X\setminus Y$ and, because $Y$ is an order ideal, these are precisely
the paths in $X$ which \emph{start} at a vertex of $X\setminus Y$. In
other words, $I_Y=Ae$. In particular, $I_Y$ is projective as a left
$A$\nobreakdash-module and \pref{p:hepi-is-idemp}~tells us that it is a
homological ideal. Since $kX/I_Y\cong kY$, the long exact sequence
of~\pref{p:pre-happel} is then
  \begin{multline}\label{eq:sing}
  \xymatrix@C-10pt{
    0 \ar[r]
      & H^0(kX,I_Y) \ar[r]
      & \HH^0(kX) \ar[r]
      & \HH^0(kY) \ar[r]
      & \Ext_{(kX)^e}^1(kX,I_Y) \ar[r]
      & \cdots
      } \\
  \xymatrix@C-17pt{
    \cdots \ar[r]
      & \Ext_{(kX)^e}^p(kX,I_Y) \ar[r]
      & \HH^p(kX) \ar[r]
      & \HH^p(kY) \ar[r]
      & \Ext_{(kX)^e}^{p+1}(kX,I_Y) \ar[r]
      & \cdots
      }
  \end{multline}
A well know result of Gerstenhaber and
Schack~\cite{GerstenhaberSchack2} states that $\HH^\*(kX)$ is
canonically isomorphic to the simplicial cohomology of the geometric
realization $\abs{X}$ of $X$ and, of course, a similar statement holds
for $kY$. Using the technique of~\cite{GerstenhaberSchack2}, one can
easily show that $\Ext_{(kX)^e}^\*(kX,I_Y) \cong H^\*(\abs{X},
\abs{Y})$, the simplicial cohomology of the pair $(\abs{X},\abs{Y})$.
Moreover, under these isomorphisms the long exact
sequence~\eqref{eq:sing} corresponds to the long exact sequence for
the cohomology of the pair~$(\abs{X},\abs{Y})$.

More generally, let $X$ be a finite poset as before and let now
$Y\subset X$ be an arbitrary subset. Let $\Ch(X)$ and $\Ch(Y)$ be the
posets of chains of $X$ and $Y$, respectively. Then $\Ch(Y)$ is an
order ideal in $\Ch(X)$. Recalling that the simplicial complex which
realizes $\Ch(X)$ is the barycentric subdivision of the one which
realizes $X$ and that simplicial cohomology is invariant under such
subdivisions, we see that $\HH^\*(k\!\Ch(X))\cong\HH^\*(kX)$. Up to
these isomorphisms, the long exact sequence~\eqref{eq:sing} provides a
long exact sequence relating $\HH^\*(kX)$, $\HH^\*(kY)$ and
  \[
  \Ext_{(k\!\Ch(X))^e}^\*(k\Ch(X), I_{\Ch(Y)}),
  \]
which again is isomorphic to a relative cohomology group. We remark
that it would be useful to have a description of these cohomology
groups directly in terms of $kX$ and $kY$.

In any case, we see that the long exact sequences for simplicial
cohomology of pairs of finite simplicial complexes are all special
cases of~\pref{p:pre-happel}.

\paragraph As observed in~\pref{p:hepi-is-idemp}, if $\phi:A\to B$ is
a homological epimorphism, the idempotency of $I=\ker\phi$ follows
from the vanishing of $\Tor^A_1(B,B)$. Looking at what happens in
degree two we find the following lemma:

\begin{Lemma*}\label{p:deg2}\cite{APT}*{Lemma~1.4}
If $\phi:A\to B$ is a homological epimorphism and $I=\ker\phi$, then
the multiplication map $I\ot_AI\to I$ is an isomorphism.~\qed
\end{Lemma*}

\medskip

\noindent Thus a homological ideal is idempotent \emph{as a
bimodule}.

\paragraph One can give various kinds of combinatorial conditions on
ideals of algebras given by quivers and relations that ensure that
they are homological. We give as a simple instance a partial converse
of~\pref{p:deg2}:

\begin{Proposition*}\label{p:ex-nocycles}
Let $Q$ be a finite quiver with path algebra $kQ$, let $J\lhd kQ$ be
an admissible ideal and put $A=kQ/J$. Let $e\in Q_0$ be a vertex such
that there are no oriented circuits in $Q$ starting in $e$. Put
$I=AeA$ and $B=A/I$. Then $I$ is a homological ideal iff
multiplication gives an isomorphism $\mu:I\ot_AI\to I$.
\end{Proposition*}

\begin{proof}
Necessity follows from~\pref{p:deg2}. In order to show the sufficiency
of the condition we need only show that $\Tor^A_+(B,B)=0$. Using twice
the long exact sequences for the functor $\Tor^A_\*$ corresponding to
the short exact sequence of $A$\nobreakdash-bimodules $0\to I\to A\to
B\to 0$, we see that
  \[
  \Tor^A_p(B,B) \cong
    \begin{cases}
      B, & \text{if $p=0$;} \\
      I/I^2, & \text{if $p=1$;} \\
      \Tor^A_1(B,I), & \text{if $p=2$;} \\
      \Tor^A_{p-2}(I,I), & \text{if $p\geq3$.}
    \end{cases}
  \]
and that there is an exact sequence
  \[\xymatrix{
  0 \ar[r]
    & \Tor^A_1(B,I) \ar[r]
    & I\ot_AI \ar[r]^-\mu
    & I \ar[r]
    & 0
  }\]
The idempotency of~$I$ and the hypothesis on~$\mu$ imply then that 
we will be done if we show that $\Tor^A_+(I,I)=0$.

Let $E$ be the subalgebra of $A$ generated by the vertices and let
$\r=\rad A$ be the Jacobson radical. Claude Cibils has shown
in~\cite{Cibils} that $A$ has a projective resolution as a
$A$\nobreakdash-bimodule of the form 
  \begin{multline*}
    \xymatrix@C-10pt{
    \cdots \ar[r]
      & A\ot_E\r^{\ot_Ep}\ot_EA \ar[r]
      & A\ot_E\r^{\ot_E(p-1)}\ot_EA \ar[r]
      & \cdots
    } \\
    \xymatrix@C-10pt{
    \cdots \ar[r]
      & A\ot_E\r^{\ot_E2}\ot_EA \ar[r]
      & A\ot_E\r\ot_EA \ar[r]
      & A\ot_E A \ar[r]
      & A \ar[r]
      & 0
    }
  \end{multline*}
It follows that $\Tor^A_\*(I,I)$ is the homology of a complex which
for each $p\geq0$ has degree~$p$ component given by
$I\ot_E\r^{\ot_Ep}\ot_EI$.

Now, since $\r$ is the ideal generated by the arrows of $Q$ and there
are no oriented circuits in $Q$ based at $e$, it is clear that 
$I\ot_E\r^{\ot_Ep}\ot_EI=0$ if $p>0$. The homology of the complex in
question is thus trivially computable and we see at once that
$\Tor^A_+(I,I)=0$.
\end{proof}

\paragraph The following proposition is
essentially~\cite{delaPenhaXi}*{Proposition~3.2.(c)} except that we do
not assume that the ideal $I$ is homological. The proof is in fact
exactly the same as the one given there, but we include it for
completeness.

\begin{Proposition*}\label{p:happel}
Let $A$ be a finite dimensional $k$\nobreakdash-algebra and let $e\in A$ be an
idempotent. Put $I=AeA$ and assume the multiplication in $A$ induces
an isomorphism of $A$\nobreakdash-bimodules $Ae\ot eA\to AeA=I$. Let
$D(\place)=\hom_k(\place, k)$ be the usual duality and let $B=A/I$. 
Then $I$ is a homological ideal and there is a long exact sequence
\begingroup
\small
  \begin{multline*}
  \xymatrix@C-5pt{
    0 \ar[r]
      & \Z(A)\cap I \ar[r]
      & \HH^0(A) \ar[r]
      & \HH^0(B) \ar[r]
      & \Ext_{A}^1(D(eA),Ae) \ar[r]
      & \cdots
      } \\
  \xymatrix@C-14.5pt{
    \cdots \ar[r]
      & \Ext_{A}^p(D(eA),Ae) \ar[r]
      & \HH^p(A) \ar[r]
      & \HH^p(B) \ar[r]
      & \Ext_{A}^{p+1}(D(eA),Ae) \ar[r]
      & \cdots
      }
  \end{multline*}
\endgroup
\end{Proposition*}

\noindent Notice that the first claim is a generalization
of~\pref{p:ex-nocycles}.

\begin{proof}
The hypothesis on $e$ implies that $I$ is idempotent and projective,
so it is homological by~\pref{p:hepi-is-idemp}. The long exact
sequence whose existence is claimed in the statement is the one
from~\pref{p:pre-happel} up to identifications. 

First of all, $H^0(A,I)\cong\hom_{A^e}(A,I)\cong\Z(A)\cap I$, so this
takes care of the beginning of the sequence. Next, dualizing the isomorphism
$I\cong Ae\ot eA$, we see that
  \begin{equation}\label{eq:di}
  D(I)\cong D(eA)\ot D(Ae),
  \end{equation}
so we have a chain of isomorphisms
  \begin{align*}
  D(\Ext_{A^e}^p(A,I))
    &  \cong \Tor^{A^e}_p(A,D(I))
    && \text{by \cite{CartanEilenberg}*{\textsc{IX}, ex.~8}} \\
    &  \cong \Tor^{A^e}_p(A,D(eA)\ot D(Ae))
    && \text{by~\eqref{eq:di}} \\
    &  \cong \Tor^A_p(D(Ae),D(eA))
    && \text{by \cite{CartanEilenberg}*{Corol.~\textsc{IX}.4.4}} \\
    &  \cong D(\Ext_A^p(D(eA),Ae))
    && \text{by \cite{CartanEilenberg}*{Prop.~\textsc{VI}.5.3}}.
  \end{align*}
Dualizing again, we see that
$\Ext_{A^e}^p(A,I)\cong\Ext_A^p(D(eA),Ae)$.
\end{proof}

\paragraph From this proposition one can easily obtain the long exact
sequence that Dieter Happel constructed in~\cite{Happel} for a
one-point extension of a finite dimensional algebra. Indeed, let $B$
be a finite dimensional $k$\nobreakdash-algebra, let $M\in\lModf B$ be
non-zero and consider the matrix algebra
  \[
  A = \begin{pmatrix} 
        B & M \\
        0 & k
      \end{pmatrix}.
  \]
If $e=\left(\begin{smallmatrix}0&0\\0&1\end{smallmatrix}\right)$, then
one easily sees that $I=AeA$ satisfies the hypothesis
of~\pref{p:happel}; moreover, it is clear that $A/I\cong B$. 
Consider now the long exact sequence corresponding to
the obvious short exact sequence
  \[\xymatrix{
  0 \ar[r]
    & M \ar[r]
    & Ae \ar[r]
    & D(eA) \ar[r]
    & 0
  }\]
and the functor $\hom_A(\place,Ae)$.

First, if $p\geq2$, we have a portion of that sequence that reads
  \begin{multline*}
    \xymatrix@C-5pt{
    \cdots \ar[r]
      & \Ext_A^{p-1}(Ae,Ae) \ar[r]
      & \Ext_A^{p-1}(M,Ae) \ar[r]
      & {}
    }
    \\
    \xymatrix@C-5pt{
    {} \ar[r]
      & \Ext_A^{p}(D(eA),Ae) \ar[r]
      & \Ext_A^{p}(Ae,Ae) \ar[r]
      & \cdots
    }
  \end{multline*}
so, since $Ae$ es projective, we conclude that
  \(
  \Ext_A^{p}(D(eA),Ae) 
  \cong
  \Ext_A^{p-1}(M,Ae)
  \).

Second, the beginning of that long exact sequence is
  \begin{multline*}
    \xymatrix@C-5pt{
    0 \ar[r]
      & \hom_A(D(eA),Ae) \ar[r]
      & \hom_A(Ae,Ae) \ar[r]
      & {}
    }
    \\
    \xymatrix@C-5pt{
    {} \ar[r]
      & \hom_A(M,Ae) \ar[r]
      & \Ext_A^1(D(eA),Ae) \ar[r]
      & 0
    }
  \end{multline*}
Using this and the fact that $\hom_A(D(eA),Ae)=0$ and $\hom_A(Ae,Ae)\cong
eAe\cong k$, we see that $\Ext_A^1(D(eA),Ae) \cong \hom_A(M,Ae) / k$.

Finally, since $\lModf B$ is a convex subcategory of $\lModf A$ and
$\hom_A(B,Ae)\cong M$, there is an isomorphism
$\Ext_A^\*(M,Ae)\cong\Ext_B^\*(M,M)$.

Using the isomorphisms thus obtained and the fact that $\Z(A)\cap
I=0$, the long exact sequence in~\pref{p:happel} becomes
\begingroup
\small
  \begin{multline*}
  \xymatrix@C-5pt{
    0 \ar[r]
      & \HH^0(A) \ar[r]
      & \HH^0(B) \ar[r]
      & \Ext_{B}^1(M,M)/k \ar[r]
      & \cdots
      } \\
  \xymatrix@C-14.5pt{
    \cdots \ar[r]
      & \Ext_{B}^p(M,M) \ar[r]
      & \HH^p(A) \ar[r]
      & \HH^p(B) \ar[r]
      & \Ext_{B}^{p+1}(M,M) \ar[r]
      & \cdots
      }
  \end{multline*}
\endgroup
which is, precisely, Happel's long exact sequence. This derivation is
explained in~\cite{delaPenhaXi}.

\paragraph The long exact sequence of Happel allows us to obtain
information on the Hochschild cohomology of an algebra which can be
built in steps from a simpler one by doing one-point extensions and
coextensions. Indeed, in terms of algebras given by quivers and
relations, one-point extensions correspond to the process of adding a
new vertex `at the top' of the quiver\dash---dually, one point
coextensions correspond to adding a new vertex `at the
bottom'\dash---and sufficiently simple quivers can be constructed
inductively starting from a vertex by adding new vertices both on the
top and on the bottom. See \cites{Gatica1, Gatica2, Gatica3} for
examples on how this inductive procedure for the computation of
cohomology is carried out.

From this point of view \pref{p:happel} becomes interesting: it allows
us to add vertices `in the middle' in certain situations, enlarging
the scope for such inductive calculations. The following proposition
is a very simple example of this.

\paragraph\label{p:ex-monomial} If $Q$ is a quiver, $e\in Q_0$ a
vertex and $\gamma=\alpha_1\cdots\alpha_k$ a path in $Q$, we say that
\definiendum{$\gamma$ has $e$ as an internal vertex} if $e$ is either
the source of one of $\alpha_1,\dots,\alpha_{k-1}$ or the target of
one of $\alpha_2,\dots,\alpha_k$. 

\begin{Proposition*}
Let $Q$ be a finite quiver with path algebra $kQ$, let $J\lhd kQ$ be
an admissible monomial ideal and put $A=kQ/J$. Let $e\in Q_0$ be a
vertex in $Q$ such that no minimal generator of $J$ has $e$ as an
internal vertex. Then $I=AeA$ is a homological ideal
and~\pref{p:happel} provides a long exact sequence relating
$\HH^\*(A)$ and $\HH^\*(A/I)$.
\end{Proposition*}

\medskip

In fact, it is not difficult to show that the condition on $e$ given
in the proposition is actually necessary for the ideal $I$ to be
homological in this case.

\begin{proof}
This follows from~\pref{p:happel} if we can prove that multiplication
induces an isomorphism $Ae\ot eA\to AeA$. The hypothesis on $e$ is
precisely what is needed for this.
\end{proof}

\paragraph As a toy example of how one can use
proposition~\pref{p:ex-monomial}, consider the algebra $A$ obtained as
the quotient of the path algebra of the following quiver by the ideal
generated by the dotted path.
  
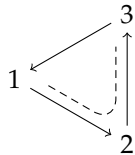
\[
    \begin{tikzpicture}
    \node (1) at (180:1) {$1$} ;
    \node (2) at (-60:1) {$2$} ;
    \node (3) at (60:1) {$3$} ;
    \draw[->] (1) -- (2) ;
    \draw[->] (2) -- (3) ;
    \draw[->] (3) -- (1) ;
    \draw[rounded corners=10pt, densely dashed, shorten >=5pt, shorten <=5pt] 
      (180:0.7) -- (-60:0.7) -- (60:0.7) ;
    \end{tikzpicture}
  \]
It is clear that the idempotent $e_1$ corresponding to the vertex $1$
satisfies the condition of the proposition, so the the ideal $I=Ae_1A$
is homological. The algebra $B=A/I$ is the path algebra of the convex
subquiver spanned by the vertices~$2$ and~$3$, which is a tree, so
$\HH^\*(B)\cong k$. Now, the $A$\nobreakdash-module $D(e_1A)$ has a
projective resolution
  \[\xymatrix{
  0 \ar[r]
    & P_3 \ar[r]
    & P_2 \ar[r]
    & P_2 \ar[r]
    & D(e_1A) \ar[r]
    & 0
  }\]
where $P_i$ is the indecomposable projective module corresponding
to the vertex~$i$, the morphism $P_2\to P_2$ maps the top onto the
socle of $P_2$ and the morphism $P_3\to P_2$ is the inclusion.
Applying the functor $\hom_A(\place,Ae_1)$ and computing, 
one easily concludes that
  \[
  \Ext_A^p(D(e_1A),Ae_1)\cong 
    \begin{cases}
    k & \text{if $p=0$ or $p=1$;} \\
    0 & \text{if $p\geq2$.}
    \end{cases}
  \]

The long exact sequence of~\pref{p:happel} then reduces in this case
to isomorphisms $\HH^p(A)=0$ for all $p\geq2$ and an exact sequence
  \begin{multline*}
    \xymatrix{
    0 \ar[r]
      & \Z(A)\cap I \ar[r]
      & \HH^0(A) \ar[r]
      & \HH^0(B) \ar[r]
      & {}
    }
      \\
    \xymatrix{
    {} \ar[r]
      & \Ext_A^1(D(e_1A),Ae_1) \ar[r]
      & \HH^1(A) \ar[r]
      & 0
    }
  \end{multline*}
The map $\HH^0(A)\to\HH^0(B)$ is surjective, so we see that
$\HH^0(A)\cong k^2$ and $\HH^1(A)\cong\Ext_A^1(D(e_1A),Ae_1)\cong k$.

Notice that $A$ is not a one-point (co)extension, so one cannot use
the classical Happel long exact sequence to compute $\HH^\*(A)$.


\vfill


\begin{bibdiv}
\begin{biblist}

\bib{APT}{article}{                                                                          
   author={Auslander, Maurice},
   author={Platzeck, Mar{\'i}a In{\'e}s},
   author={Todorov, Gordana},
   title={Homological theory of idempotent ideals},
   journal={Trans. Amer. Math. Soc.},
   volume={332},
   date={1992},
   number={2},
   pages={667--692},
   issn={0002-9947},
}

\bib{Bardzell}{article}{
   author={Bardzell, Michael J.},
   title={The alternating syzygy behavior of monomial algebras},
   journal={J. Algebra},
   volume={188},
   date={1997},
   number={1},
   pages={69--89},
   issn={0021-8693},
}

\bib{BLM}{article}{
   author={Bardzell, Michael J.},
   author={Locateli, Ana Claudia},
   author={Marcos, Eduardo N.},
   title={On the Hochschild cohomology of truncated cycle algebras},
   journal={Comm. Algebra},
   volume={28},
   date={2000},
   number={3},
   pages={1615--1639},
   issn={0092-7872},
}
\bib{BACH}{article}{
   author={Buenos Aires Cyclic Homology Group},
   title={Cyclic homology of algebras with one generator},
   journal={$K$\nobreakdash-Theory},
   volume={5},
   date={1991},
   number={1},
   pages={51\ndash 69},
   issn={0920-3036},
}

\bib{CartanEilenberg}{book}{
   author={Cartan, Henri},
   author={Eilenberg, Samuel},
   title={Homological algebra},
   publisher={Princeton University Press},
   place={Princeton, N. J.},
   date={1956},
   pages={xv+390},
}

\bib{Cibils1}{article}{
   author={Cibils, Claude},
   title={Tensor Hochschild homology and cohomology},
   conference={
      title={Interactions between ring theory and representations of
      algebras (Murcia)},
   },
   book={
      series={Lecture Notes in Pure and Appl. Math.},
      volume={210},
      publisher={Dekker},
      place={New York},
   },
   date={2000},
   pages={35--51},
}

\bib{Cibils}{article}{
   author={Cibils, Claude},
   title={Tensor Hochschild homology and cohomology},
   conference={
      title={Interactions between ring theory and representations of
      algebras (Murcia)},
   },
   book={
      series={Lecture Notes in Pure and Appl. Math.},
      volume={210},
      publisher={Dekker},
      place={New York},
   },
   date={2000},
   pages={35--51},
}

\bib{Gatica1}{article}{
   author={Gatica, Mar{\'{\i}}a Andrea},
   author={Rey, Andrea Alejandra},
   title={Computing the Hochschild cohomology groups of some families of
   incidence algebras},
   journal={Comm. Algebra},
   volume={34},
   date={2006},
   number={6},
   pages={2039--2056},
   issn={0092-7872},
}

\bib{Gatica2}{article}{
   author={Gatica, Mar{\'{\i}}a Andrea},
   author={Redondo, Mar{\'{\i}}a Julia},
   title={Hochschild cohomology of incidence algebras as one-point
   extensions},
   note={Special issue on linear algebra methods in representation theory},
   journal={Linear Algebra Appl.},
   volume={365},
   date={2003},
   pages={169--181},
   issn={0024-3795},
}

\bib{Gatica3}{article}{
   author={Gatica, Mar{\'{\i}}a Andrea},
   author={Redondo, Mar{\'{\i}}a Julia},
   title={Hochschild cohomology and fundamental groups of incidence
   algebras},
   journal={Comm. Algebra},
   volume={29},
   date={2001},
   number={5},
   pages={2269--2283},
   issn={0092-7872},
}

\bib{GeigleLenzig}{article}{
   author={Geigle, Werner},
   author={Lenzing, Helmut},
   title={Perpendicular categories with applications to representations
   and sheaves},
   journal={J. Algebra},
   volume={144},
   date={1991},
   number={2},
   pages={273\ndash 343},
   issn={0021-8693},
}

\bib{Gerstenhaber}{article}{
   author={Gerstenhaber, Murray},
   title={The cohomology structure of an associative ring},
   journal={Ann. of Math. (2)},
   volume={78},
   date={1963},
   pages={267\ndash 288},
   issn={0003-486X},
}

\bib{GerstenhaberSchack1}{article}{
   author={Gerstenhaber, Murray},
   author={Schack, Samuel D.},
   title={Algebras, bialgebras, quantum groups, and algebraic deformations},
   booktitle={Deformation theory and quantum groups with applications to
        mathematical physics (Amherst, MA, 1990)},
   series={Contemp. Math.},
   volume={134},
   pages={51\ndash 92},
   publisher={Amer. Math. Soc.},
   place={Providence, RI},
   date={1992},
}

\bib{GerstenhaberSchack2}{article}{
   author={Gerstenhaber, Murray},
   author={Schack, Samuel D.},
   title={Simplicial cohomology is Hochschild cohomology},
   journal={J. Pure Appl. Algebra},
   volume={30},
   date={1983},
   number={2},
   pages={143--156},
   issn={0022-4049},
}

\bib{GreenSolberg}{article}{
   author={Green, Edward L.},
   author={Solberg, {\oldO}yvind},
   title={Hochschild cohomology rings and triangular rings},
   conference={
      title={Representations of algebra. Vol. I, II},
   },
   book={
      publisher={Beijing Norm. Univ. Press, Beijing},
   },
   date={2002},
   pages={192--200},
}

\bib{GreenMarcosSnashall}{article}{
   author={Green, Edward L.},
   author={Marcos, Eduardo N.},
   author={Snashall, Nicole},
   title={The Hochschild cohomology ring of a one point extension},
   journal={Comm. Algebra},
   volume={31},
   date={2003},
   number={1},
   pages={357--379},
   issn={0092-7872},
}

\bib{Happel}{article}{
   author={Happel, Dieter},
   title={Hochschild cohomology of finite-dimensional algebras},
   booktitle={S\'eminaire d'Alg\`ebre Paul Dubreil et Marie-Paul Malliavin,
        39\`eme Ann\'ee (Paris, 1987/1988)},
   series={Lecture Notes in Math.},
   volume={1404},
   pages={108\ndash 126},
   publisher={Springer},
   place={Berlin},
   date={1989},
}

\bib{HiltonRees}{article}{
   author={Hilton, Peter John},
   author={Rees, David},
   title={Natural maps of extension functors and a theorem of R. G. Swan},
   journal={Proc. Cambridge Philos. Soc.},
   volume={57},
   date={1961},
   pages={489\ndash 502},
}

\bib{HiltonStammbach}{book}{
   author={Hilton, Peter John},
   author={Stammbach, Urs},
   title={A course in homological algebra},
   series={Graduate Texts in Mathematics},
   volume={4},
   publisher={Springer-Verlag},
   place={New York},
   date={1997},
   pages={xii+364},
   isbn={0-387-94823-6},
}

\bib{Hochschild}{article}{
   author={Hochschild, Gerhard},
   title={On the cohomology groups of an associative algebra},
   journal={Ann. of Math. (2)},
   volume={46},
   date={1945},
   pages={58\ndash 67},
   issn={0003-486X},
}

\bib{Holm:Brauer}{article}{
   author={Holm, Thorsten},
   title={Hochschild cohomology of Brauer tree algebras},
   journal={Comm. Algebra},
   volume={26},
   date={1998},
   number={11},
   pages={3625--3646},
   issn={0092-7872},
}
	
\bib{Holm}{article}{
   author={Holm, Thorsten},
   title={Hochschild cohomology rings of algebras $k[X]/(f)$},
   journal={Beitr\"age Algebra Geom.},
   volume={41},
   date={2000},
   number={1},
   pages={291\ndash 301},
   issn={0138-4821},
}

\bib{Karoubi}{article}{
   author={Karoubi, Max},
   title={Homologie cyclique et $K$-th\'eorie},
   journal={Ast\'erisque},
   number={149},
   date={1987},
   pages={147},
   issn={0303-1179},
}

\bib{MacLane}{book}{
   author={MacLane, Saunders},
   title={Homology},
   publisher={Springer-Verlag},
   place={Berlin},
   date={1967},
   pages={x+422},
}

\bib{MichelenaPlatzeck}{article}{
   author={Michelena, Sandra},
   author={Platzeck, Mar{\'{\i}}a In{\'e}s},
   title={Hochschild cohomology of triangular matrix algebras},
   journal={J. Algebra},
   volume={233},
   date={2000},
   number={2},
   pages={502--525},
   issn={0021-8693},
}

\bib{delaPenhaXi}{article}{
   author={de la Pe{\~n}a, Jos{\'e} Antonio},
   author={Xi, Changchang},
   title={Hochschild cohomology of algebras with homological ideals},
   journal={Tsukuba J. Math.},
   volume={30},
   date={2006},
   number={1},
   pages={61--79},
   issn={0387-4982},
}

\bib{Quillen}{article}{
   author={Quillen, Daniel},
   title={The spectrum of an equivariant cohomology ring. I, II},
   journal={Ann. of Math. (2)},
   volume={94},
   date={1971},
   pages={549--572; ibid. (2) 94 (1971), 573--602},
   issn={0003-486X},
   review={\MR{0298694 (45 \#7743)}},
}

\bib{Stanley}{book}{
   author={Stanley, Richard Peter},
   title={Enumerative combinatorics. Vol. 1},
   series={Cambridge Studies in Advanced Mathematics},
   volume={49},
   publisher={Cambridge University Press},
   place={Cambridge},
   date={1997},
   pages={xii+325},
   isbn={0-521-55309-1},
   isbn={0-521-66351-2},
}

\end{biblist}
\end{bibdiv}
\end{document}